\documentclass[12pt]{article}
\usepackage[english]{babel}
\usepackage[T1]{fontenc}
\usepackage[utf8]{inputenc}
\usepackage[top = 1in, bottom = 1in, left = 1in, right = 1in]{geometry}
\usepackage{amsmath,amsxtra,amssymb,latexsym,amscd,amsthm,enumitem,calligra,mathrsfs}
\usepackage{graphicx}
\usepackage[mathscr]{eucal}
\usepackage{fancybox}
\usepackage[all]{xy}
\usepackage{tikz}
\usepackage{tikz-cd}
\usepackage{amscd}
\usepackage{mathtools}
\usepackage{stmaryrd}
\usepackage[backend=biber,style=alphabetic]{biblatex}
\usetikzlibrary{positioning}
\usetikzlibrary{arrows}
\usepackage{hyperref}
\usepackage[nottoc]{tocbibind}

\usepackage{imakeidx}
\makeindex[columns=2, title=Index, intoc]

\DeclareFontFamily{U}{wncy}{}
\DeclareFontShape{U}{wncy}{m}{n}{<->wncyr10}{}
\DeclareSymbolFont{mcy}{U}{wncy}{m}{n}
\DeclareMathSymbol{\Sh}{\mathord}{mcy}{"58}

\newcommand{\Rbb}{\mathbb{R}}
\newcommand{\Cbb}{\mathbb{C}}

\newcommand{\Qbb}{\mathbb{Q}}
\newcommand{\Zbb}{\mathbb{Z}}
\newcommand{\Fpbb}{\mathbb{F}_{p}}

\newcommand{\Qpbb}{\mathbb{Q}_{p}}
\newcommand{\Zpbb}{\mathbb{Z}_{p}}
\newcommand{\Pbb}{\mathbb{P}}
\newcommand{\Abb}{\mathbb{A}}
\newcommand{\Gmbb}{\mathbb{G}_{m}}

\newcommand{\ds}{\displaystyle}

\theoremstyle{definition}
\newtheorem{definition}{Definition}[section]
\newtheorem{remark}[definition]{Remark}
\newtheorem{example}[definition]{Example}

\newtheorem{assumption}[definition]{Assumption}

\theoremstyle{plain}
\newtheorem{theorem}{Theorem}[section]
\newtheorem{proposition}[theorem]{Proposition}
\newtheorem{lemma}[theorem]{Lemma}
\newtheorem{corollary}[theorem]{Corollary}
\newtheorem{conjecture}[theorem]{Conjecture}

\begin{document}




\title{\textbf{Brauer--Manin obstruction for integral points on Markoff-type cubic surfaces}}
	\author{\textsc{Quang-Duc DAO}} 
	\date{}
	\maketitle

\begin{abstract}
Following \cite{GS22}, \cite{LM20} and \cite{CTWX20}, we study the Brauer--Manin obstruction for integral points on similar Markoff-type cubic surfaces. In particular, we construct a family of counterexamples to strong approximation which can be explained by the Brauer--Manin obstruction with some counting results of similar nature to those in \cite{LM20} and \cite{CTWX20}. We also give some counterexamples to the integral Hasse principle which cannot be explained by the (algebraic) Brauer--Manin obstruction.
\end{abstract}



\section{Introduction}
Let $X$ be an affine variety over $\Qbb$, and $\mathcal{X}$ an integral model of $X$ over $\Zbb$, i.e. an affine scheme of finite type over $\Zbb$ whose generic fiber is isomorphic to $U$. Define the set of adelic points $X(\textbf{\textup{A}}_{\Qbb}) := \sideset{}{'}\prod_{p} X(\Qpbb)$, where $p$ is a prime number or $p = \infty$ (with $\Qbb_{\infty} = \Rbb$). Similarly, define $\mathcal{X}(\textbf{\textup{A}}_{\Zbb}) := \prod_{p} \mathcal{X}(\Zpbb)$ (with $\Zbb_{\infty} = \Rbb$). We say that $X$ \textit{fails the Hasse principle} if 
$$ X(\textbf{\textup{A}}_{\Qbb}) \not= \emptyset \hspace{0.5cm}\textup{but}\hspace{0.5cm} X(\Qbb) = \emptyset $$
We say that $\mathcal{X}$ \textit{fails the integral Hasse principle} if 
$$ \mathcal{X}(\textbf{\textup{A}}_{\Zbb}) \not= \emptyset \hspace{0.5cm}\textup{but}\hspace{0.5cm} \mathcal{X}(\Zbb) = \emptyset. $$
We say that $X$ \textit{satisfies weak approximation} if the image of $X(\Qbb)$ in $\prod_{v} X(\Qbb_{v})$ is dense, where the product is taken over all places of $\Qbb$. Finally, we say that $\mathcal{X}$ \textit{satisfies strong approximation} if $\mathcal{X}(\Zbb)$ is dense in $\mathcal{X}(\textbf{\textup{A}}_{\Zbb})_{\bullet} := \prod_{p} \mathcal{X}(\Zpbb) \times \pi_{0}(X(\Rbb))$, where $\pi_{0}(X(\Rbb))$ denotes the set of connected components of $X(\Rbb)$. Note that we work with $\pi_{0}(X(\Rbb))$ since $\mathcal{X}(\Zbb)$ is never dense in $X(\Rbb)$ for topological reasons (see \cite[Example 2.2]{Con12}).

In general, the Hasse principle for varieties does not hold. In his 1970 ICM address \cite{Man71}, Manin introduced a natural cohomological obstruction to the Hasse principle, namely the \textbf{Brauer--Manin obstruction} (which has been extended to its integral version in \cite{CTX09}). If $\textup{Br}\,X$ denotes the cohomological Brauer group of $X$, i.e. $\textup{Br}\,X := \textup{H}^{2}_{\textup{ét}}(X,\Gmbb)$, we have a natural pairing from class field theory:
$$ X(\textbf{\textup{A}}_{\Qbb}) \times \textup{Br}\,X \rightarrow \Qbb/\Zbb. $$
If we define $X(\textbf{\textup{A}}_{\Qbb})^{\textup{Br}}$ to be the left kernel of this pairing, then the exact sequence of Albert--Brauer--Hasse--Noether gives us the relation:
$$ X(\Qbb) \subseteq X(\textbf{\textup{A}}_{\Qbb})^{\textup{Br}} \subseteq X(\textbf{\textup{A}}_{\Qbb}). $$
Similarly, by defining the Brauer--Manin set $\mathcal{X}(\textbf{\textup{A}}_{\Zbb})_{\bullet}^{\textup{Br}}$, we also have that 
$$ \mathcal{X}(\Zbb) \subseteq \mathcal{X}(\textbf{\textup{A}}_{\Zbb})_{\bullet}^{\textup{Br}} \subseteq \mathcal{X}(\textbf{\textup{A}}_{\Zbb})_{\bullet}. $$
This gives the so-called \textit{integral} Brauer--Manin obstruction. We say that \textit{the Brauer--Manin obstruction to the (resp. integral) Hasse principle is the only one} if 

$$ X(\textbf{\textup{A}}_{\Qbb})^{\textup{Br}} \not= \emptyset \iff X(\Qbb) \not= \emptyset.$$

$$ (\mathcal{X}(\textbf{\textup{A}}_{\Zbb})_{\bullet}^{\textup{Br}} \not= \emptyset \iff \mathcal{X}(\Zbb) \not= \emptyset.)$$
If there is no confusion, we can omit the symbol $\bullet$ for the set of local integral points and the corresponding Brauer--Manin set.
\~\\
\indent We are particularly interested in the case where $X$ is a hypersurface, defined by a polynomial equation of degree $d$ in an affine space. The case $d = 1$ is easy. The case $d = 2$ considers the arithmetic of quadratic forms: for rational points, the Hasse principle is always satisfied by the Hasse--Minkowski theorem, and for integral points, the Brauer--Manin obstruction to the integral Hasse principle is the only one (up to an isotropy assumption) due to work of Colliot-Thélène, Xu \cite{CTX09} and Harari \cite{Ha08}. However, the case $d = 3$ (of cubic hypersurfaces) is still largely open, especially for integral points. Overall, the arithmetic of integral points on the affine cubic surfaces over number fields is still little understood. For example, the question to determine which integers can be written as sums of three cubes of integers is still open. In this first problem, for the affine variety defined by the equation 
$$ x^{3} + y^{3} + z^{3} = a, $$
where $a$ is a fixed integer, Colliot-Thélène and Wittenberg in \cite{CTW12} proved that there is no Brauer--Manin obstruction to the integral Hasse principle (if $a$ is not of the form $9n \pm 4$). However, the existence of such an integer $a$ remains unknown in general, with the first example now being the case when $a = 114$. On the other hand, in a related problem, there is no Brauer--Manin obstruction to the existence of an integral point on the cubic surface defined by 
$$ x^{3} + y^{3} + 2z^{3} = a, $$ for any $a \in \Zbb$, also proven in \cite{CTW12}.

Another interesting example of affine cubic surfaces that we consider is given by \textbf{Markoff surfaces} $U_{m}$ which are defined by 
$$ x^{2} + y^{2} + z^{2} - xyz = m, $$ where $m$ is an integer parameter. The very first (original) class of Markoff surfaces which was studied is the one given by this equation with $m = 0$ in a series of papers \cite{BGS16a}, \cite{BGS16b}, and recently \cite{Che21}. They study \emph{strong approximation mod $p$} for $U_{0}$ for any prime $p$ and present a Strong Approximation Conjecture (\cite[Conjecture 1]{BGS16b}):
\begin{conjecture}
For any prime $p$, $U_{0}(\Zbb/p\Zbb)$ consists of two $\Gamma$ orbits, namely $\{(0, 0, 0)\}$ and $U_{0}^{*}(\Zbb/p\Zbb) = U_{0}(\Zbb/p\Zbb) \backslash \{(0, 0, 0)\}$. Here $\Gamma$ is a group
of affine integral morphisms of $\Abb^{3}$ generated by the permutations of the coordinates and the Vieta involutions.
\end{conjecture}
Combining the above three papers by Bourgain, Gamburd, Sarnak, and Chen, the Conjecture is established for all but finitely many primes (see \cite[Theorem 5.5.5]{Che21}), which also implies that $U_{0}$ satisfies strong approximation mod $p$ for all by finitely many primes.

On the other hand, in \cite{GS22}, Ghosh and Sarnak study the integral points on those affine Markoff surfaces $U_{m}$ with general $m$, both from a theoretical point of view and by numerical evidence. They prove that for almost all $m$, the integral Hasse principle holds, and that there are infinitely many $m$'s for which it fails (Hasse failures). Furthermore, their numerical experiments suggest particularly a proportion of integers $m$ satisfying $|m| \leq M$ of the power $M^{0,8875\dots+o(1)}$ for which the integral Hasse principle is not satisfied.

Subsequently, Loughran and Mitankin \cite{LM20} proved that asymptotically only a proportion of $M^{1/2}/(\log M)^{1/2}$ of integers $m$ such that $-M \leq m \leq M$ presents an integral Brauer--Manin obstruction to the Hasse principle. They also obtained a lower bound, asymptotically $M^{1/2}/\log M$, for the number of Hasse failures which cannot be explained by the Brauer--Manin obstruction. After Colliot-Thélène, Wei, and Xu \cite{CTWX20} obtained a slightly stronger lower bound than the one given in \cite{LM20}, no better result than their number $M^{1/2}/(\log M)^{1/2}$ has been known until now. In other words, with all the current results, one does not have a satisfying comparison between the numbers of Hasse failures which can be explained by the Brauer--Manin obstruction and which cannot be explained by this obstruction. Meanwhile, for strong approximation, it has been proven to almost never hold for Markoff surfaces in \cite{LM20} and then absolutely never be the case in \cite{CTWX20}. 
Here we recall an important conjecture given by Ghosh and Sarnak.

\begin{conjecture} \textup{\cite[Conjecture 10.2]{GS22}} 
The number of Hasse failures satisfies that
$$ \# \{m \in \Zbb: 0 \leq m \leq M,\ \mathcal{U}_{m}(\Abb_{\Zbb}) \not= \emptyset\ \text{but}\ \mathcal{U}_{m}(\Zbb) = \emptyset \} \approx C_{0}M^{\theta}, $$
for some $C_{0} > 0$ and some $\frac{1}{2} < \theta < 1$.
\end{conjecture}

The above conjecture also means that \emph{almost all} counterexamples to the integral Hasse principle for Markoff surfaces cannot be explained by the Brauer--Manin obstruction, thanks to the result obtained by \cite{LM20}. While the question of counting all counterexamples to the integral Hasse principle for Markoff surfaces remains largely open, we will focus on another family in this paper. More precisely, we are going to study the set of integral points of a different Markoff-type cubic surfaces whose origin is similar to that of the original Markoff surfaces $U_{m}$, namely the \textit{relative character varieties} which will be introduced in Section 2, using the Brauer--Manin obstruction as well. The surfaces are given by the cubic equation:

$$ x^{2} + y^{2} + z^{2} + xyz = ax + by + cz + d, $$

\noindent where $a,b,c,d \in \Zbb$ are parameters that satisfy some specific relations to be discussed later. Due to the similar appearance as that of the original Markoff surfaces, one may expect to find some similarities in their arithmetic as well. One of the main results in our paper is the following, saying that a positive proportion of these relative character varieties have no (algebraic) Brauer--Manin obstruction to the integral Hasse principle as well as fail \emph{strong approximation}, and those failures can be explained by the Brauer--Manin obstruction.

\begin{theorem}
Let $\mathcal{U}$ be the affine scheme over $\Zbb$ defined by 
\begin{equation*}
x^{2} + y^{2} + z^{2} + xyz = ax + by + cz + d,
\end{equation*}
where 
\begin{equation*}
	\begin{cases*}
	a = k_{1}k_{2} + k_{3}k_{4} \\
	b = k_{1}k_{4} + k_{2}k_{3} \\
	c = k_{1}k_{3} + k_{2}k_{4}
	\end{cases*}
	\hspace{0.5cm} \textup{and} \hspace{0.5cm} d = 4 - \sum_{i=1}^{4} k_{i}^{2} - \prod_{i=1}^{4} k_{i}, 
\end{equation*}
such that the projective closure $X \subset \Pbb^{3}_{\Qbb}$ of $U = \mathcal{U} \times_{\Zbb} \Qbb$ is smooth. Then we have
$$ \#\{ k = (k_{1},k_{2},k_{3},k_{4}) \in \Zbb^{4}, |k_{i}| \leq M\; \forall\; 1 \leq i \leq 4 : \emptyset \not= \mathcal{U}(\textbf{\textup{A}}_{\Zbb})^{\textup{Br}_{1}} \not= \mathcal{U}(\textbf{\textup{A}}_{\Zbb}) \} \asymp M^{4} $$
as $M \rightarrow +\infty$.
\end{theorem}

The structure of the paper is as follows. In Section 2, we provide some background on character varieties and a natural origin of the Markoff-type cubic surfaces. In Section 3, we first review some important results on affine cubic surfaces given by the complement of three coplanar lines and their Brauer groups. After the general setting, we turn our attention to the natural smooth projective compactifications of the Markoff-type cubic surfaces, where we explicitly calculate the (algebraic) Brauer group of the compactification, and then we complete the analysis of the Brauer group by calculating the Brauer group of the affine surfaces. In Section 4, we use the Brauer group to give explicit examples of Brauer--Manin obstructions to the integral Hasse principle, and give some counting results for the frequency of the obstructions. Finally, in Section 5, we make some important remarks to compare some results in this paper to those of Markoff surfaces in recent work that we follow, and We also give some counterexamples to the integral Hasse principle which cannot be explained by the Brauer--Manin obstruction.
\\~\\
\indent \textbf{Notation.} Let $k$ be a field and $\overline{k}$ a separable closure of $k$. We let $G_{k} := \textup{Gal}(\overline{k}/k)$ be the absolute Galois group. A $k$-variety is a separated $k$-scheme of finite type. If $X$ is a $k$-variety, we write $\overline{X} = X \times_{k} \overline{k}$. Let $k[X] = \textup{H}^{0}(X,\mathcal{O}_{X})$ and $\overline{k}[X] = \textup{H}^{0}(\overline{X},\mathcal{O}_{\overline{X}})$. If $X$ is an integral $k$-variety, let $k(X)$ denote the function field of $X$. If $X$ is a geometrically integral $k$-variety, let $\overline{k}(X)$ denote the function field of $\overline{X}$. 

Let $\textup{Pic}\,X = \textup{H}^{1}_{\textup{Zar}}(X,\Gmbb) = \textup{H}^{1}_{\textup{ét}}(X,\Gmbb)$ denote the Picard group of a scheme $X$. Let $\textup{Br}\,X = \textup{H}^{2}_{\textup{ét}}(X,\Gmbb)$ denote the Brauer group of $X$. Let
$$ \textup{Br}_{1}\,X := \textup{Ker}[\textup{Br}\,X \rightarrow \textup{Br}\,\overline{X}] $$
denote the \textbf{algebraic Brauer group} of a $k$-variety $X$ and let $\textup{Br}_{0}\,X \subset \textup{Br}_{1}\,X$ denote the image of $\textup{Br}\,k \rightarrow \textup{Br}\,X$. The image of $\textup{Br}\,X \rightarrow \textup{Br}\,\overline{X}$ is called the \textbf{transcendental Brauer group} of $X$.

Given a field $F$ of characteristic zero containing a primitive $n$-th root of unity $\zeta = \zeta_{n}$, we have $\textup{H}^{2}(F,\mu_{n}^{\otimes 2}) = \textup{H}^{2}(F,\mu_{n}) \otimes \mu_{n}$. The choice of $\zeta_{n}$ then defines an isomorphism $\textup{Br}(F)[n] = \textup{H}^{2}(F,\mu_{n}) \cong \textup{H}^{2}(F,\mu_{n}^{\otimes 2})$. Given two elements $f, g \in F^{\times}$, we have their classes $(f)$ and $(g)$ in $F^{\times}/F^{\times n} = \textup{H}^{1}(F,\mu_{n})$. We denote by $(f,g)_{\zeta} \in \textup{Br}(F)[n] = \textup{H}^{2}(F,\mu_{n})$ the class corresponding to the cup-product $(f) \cup (g) \in \textup{H}^{2}(F,\mu_{n}^{\otimes 2})$. Suppose $F/E$ is a finite Galois extension with Galois group $G$. Given $\sigma \in G$ and $f,g \in F^{\times}$, we have $\sigma((f,g)_{\zeta_{n}}) = (\sigma(f),\sigma(g))_{\sigma(\zeta_{n})} \in \textup{Br}(F)$. In particular, if $\zeta_{n} \in E$, then $\sigma((f,g)_{\zeta_{n}}) = (\sigma(f),\sigma(g))_{\zeta_{n}}$. For all the details, see \cite[Sections 4.6, 4.7]{GS17}.

Let $R$ be a discrete valuation ring with fraction field $F$ and residue field $\kappa$. Let $v$ denote the valuation $F^{\times} \rightarrow \Zbb$. Let $n > 1$ be an integer invertible in $R$. Assume that $F$ contains a primitive $n$-th root of unity $\zeta$. For $f,g \in F^{\times}$, we have the residue map
$$ \partial_{R} : \textup{H}^{2}(F,\mu_{n}) \rightarrow \textup{H}^{1}(\kappa,\Zbb/n\Zbb) \cong \textup{H}^{1}(\kappa,\mu_{n}) = \kappa^{\times}/\kappa^{\times n}, $$
where $\textup{H}^{1}(\kappa,\Zbb/n\Zbb) \cong \textup{H}^{1}(\kappa,\mu_{n})$ is induced by the isomorphism $\Zbb/n\Zbb \simeq \mu_{n}$ sending $1$ to $\zeta$. This map sends the class of $(f,g)_{\zeta} \in \textup{Br}(F)[n] = \textup{H}^{2}(F,\mu_{n})$ to 
$$ (-1)^{v(f)v(g)} \textup{class}(g^{v(f)}/f^{v(g)}) \in \kappa/\kappa^{\times n}. $$

For a proof of these facts, see \cite{GS17}. Here we recall some precise references. Residues in Galois cohomology with finite coefficients are defined in \cite[Construction 6.8.5]{GS17}. Comparison of residues in Milnor K-Theory and Galois cohomology is given in \cite[Proposition 7.5.1]{GS17}. The explicit formula for the residue in Milnor’s group K2 of a discretely valued field is given in \cite[Example 7.1.5]{GS17}.
\\~\\
\indent \textbf{Acknowledgements.} I thank Cyril Demarche for his help and supervision during my PhD study at the Institute of Mathematics of Jussieu. I thank Kevin Destagnol for his help with the computations in Section 4.3 using analytic number theory. I thank Vladimir Mitankin for his useful remarks and suggestions, especially regarding Section 5. I thank Jean-Louis Colliot-Thélène, Fei Xu, and Daniel Loughran for their helpful comments and encouragement. This project has received funding from the European Union’s Horizon 2020 Research and Innovation Programme under the Marie Skłodowska-Curie grant agreement No. 754362 and from the Vietnam Academy of Science and Technology's 2022 Support Programme for Junior Researchers. I thank the reviewers for their careful reading of my paper and their many insightful comments and suggestions which helped me improve considerably the manuscript.

\section{Background}
The main reference to look up notations that we use here is \cite{Wha20}, mostly Chapter 2.
\subsection{Character varieties}
First, we introduce an important origin of the Markoff-type cubic surfaces which comes from \emph{character varieties}, as studied in \cite{Wha20}. Throughout this section, an \textit{algebraic variety} is a scheme of finite type over a field. Given an affine variety $X$ over a field $k$, we denote by $k[X]$ its coordinate ring over $k$. If moreover $X$ is integral, then $k(X)$ denotes its function field over $k$. Given a commutative ring $A$ with unity, the elements of $A$ will be referred to as regular functions on the affine scheme $\textup{Spec}\,A$.

\begin{definition}
Let $\pi$ be a finitely generated group. Its ($\textup{SL}_{2}$) \textit{representation variety} $\textup{Rep}(\pi)$ is the affine scheme defined by the functor
$$ A \mapsto \textup{Hom}(\pi,\textup{SL}_{2}(A)) $$
for every commutative ring $A$. Assume that $\pi$ has a sequence of generators of $m$ elements, then we have a presentation of $\textup{Rep}(\pi)$ as a closed subscheme of $\textup{SL}_{2}^{m}$ defined by equations coming from relations among the generators. For each $a \in \pi$, let $\textup{tr}_{a}$ be the regular function on $\textup{Rep}(\pi)$ given by $\rho \mapsto \textup{tr}\,\rho(a)$.

The ($\textup{SL}_{2}$) \textit{character variety} of $\pi$ over $\Cbb$ is then defined to be the affine invariant theoretic quotient
$$ X(\pi) := \textup{Rep}(\pi) \sslash \textup{SL}_{2} = \textup{Spec}\left(\Cbb[\textup{Rep}(\pi)]^{\textup{SL}_{2}(\Cbb)}\right) $$
under the conjugation action of $\textup{SL}_{2}$. 
\end{definition}

The regular function $\textup{tr}_{a}$ for each $a \in \pi$ clearly descends to a regular function on $X(\pi)$. Furthermore, from the fact that $\textup{tr}(I_{2}) = 2$ and $\textup{tr}(A)\textup{tr}(B) = \textup{tr}(AB) + \textup{tr}(AB^{-1})$, for $I_{2} \in \textup{SL}_{2}(\Cbb)$ being the identity matrix and for any $A,B \in \textup{SL}_{2}(\Cbb)$, we can deduce a natural model of $X(\pi)$ over $\Zbb$ as the spectrum of
$$ R(\pi) := \Zbb[\textup{tr}_{a} : a \in \pi] / (\textup{tr}_{1}-2, \textup{tr}_{a}\textup{tr}_{b} - \textup{tr}_{ab} - \textup{tr}_{ab^{-1}}). $$
Given any integral domain $A$ with fraction field $F$ of characteristic zero, the $A$-points of $X(\pi)$ parametrize the Jordan equivalence classes of $\textup{SL}_{2}(F)$-representations of $\pi$ having character valued in $A$.

\begin{example}
Denote by $F_{m}$ the free group on $m \geq 1$ generators $a_{1},\dots,a_{m}$. By Goldman's results used in \cite{Wha20}, we have the following important examples:

\begin{enumerate}
	\item[(1)] $\textup{tr}_{a_{1}} : X(F_{1}) \simeq \Abb^{1}$.
	
	\item[(2)] $(\textup{tr}_{a_{1}},\textup{tr}_{a_{2}},\textup{tr}_{a_{3}}) : X(F_{2}) \simeq \Abb^{3}$.
	
	\item[(3)] The coordinate ring $\Qbb[X(F_{3})]$ is the quotient of the polynomial ring
	$$ \Qbb[\textup{tr}_{a_{1}}, \textup{tr}_{a_{2}}, \textup{tr}_{a_{3}}, \textup{tr}_{a_{1}a_{2}}, \textup{tr}_{a_{2}a_{3}}, \textup{tr}_{a_{1}a_{3}}, \textup{tr}_{a_{1}a_{2}a_{3}}, \textup{tr}_{a_{1}a_{3}a_{2}}] $$
	by the ideal generated by two elements 
	$$ \textup{tr}_{a_{1}a_{2}a_{3}} + \textup{tr}_{a_{1}a_{3}a_{2}} - (\textup{tr}_{a_{1}a_{2}} \textup{tr}_{a_{3}} + \textup{tr}_{a_{1}a_{3}} \textup{tr}_{a_{2}} + \textup{tr}_{a_{2}a_{3}} \textup{tr}_{a_{1}} - \textup{tr}_{a_{1}} \textup{tr}_{a_{2}} \textup{tr}_{a_{3}}) $$
	and
	\[
	\begin{aligned}
	\textup{tr}_{a_{1}a_{2}a_{3}} \textup{tr}_{a_{1}a_{3}a_{2}} - \{&(\textup{tr}_{a_{1}}^{2} + \textup{tr}_{a_{2}}^{2} + \textup{tr}_{a_{3}}^{2}) + (\textup{tr}_{a_{1}a_{2}}^{2} + \textup{tr}_{a_{2}a_{3}}^{2} + \textup{tr}_{a_{1}a_{3}}^{2})\\ &- (\textup{tr}_{a_{1}} \textup{tr}_{a_{2}} \textup{tr}_{a_{1}a_{2}} + \textup{tr}_{a_{2}} \textup{tr}_{a_{3}} \textup{tr}_{a_{2}a_{3}} + \textup{tr}_{a_{1}} \textup{tr}_{a_{3}} \textup{tr}_{a_{1}a_{3}})\\ &+ \textup{tr}_{a_{1}a_{2}} \textup{tr}_{a_{2}a_{3}} \textup{tr}_{a_{1}a_{3}} - 4\}.
	\end{aligned}
	\]
\end{enumerate}
\end{example}

Now given a connected smooth compact manifold $M$, we consider the \textit{moduli of local systems} on $M$ which is the character variety $X(M) := X(\pi_{1}(M))$ of its fundamental group. More generally, given a smooth manifold $M = M_{1} \sqcup \cdots \sqcup M_{m}$ with finitely many connected components $M_{i}$ for $1 \leq i \leq m$, define 
$$ X(M) := X(M_{1}) \times \cdots \times X(M_{m}). $$
The construction of the moduli space $X(M)$ is functorial in the manifold $M$. More precisely, any smooth map $f : M \rightarrow N$ of manifolds induces a morphism $f^{*} : X(N) \rightarrow X(M)$, depending only on the homotopy class of $f$, given by pullback of local systems.

Let $\Sigma$ be a surface. For each curve $a \in \Sigma$, there is a well-defined regular function $\textup{tr}_{a} : X(\Sigma) \rightarrow X(a) \simeq \Abb^{1}$, which agrees with $\textup{tr}_{\alpha}$ for any $\alpha \in \pi_{1}(\Sigma)$ freely homotopic to a parametrization of $a$. The boundary curves $\partial\Sigma$ of $\Sigma$ induce a natural morphism
$$ \textup{tr}_{\partial\Sigma} = (-)|_{\partial\Sigma} : X(\Sigma) \rightarrow X(\partial\Sigma). $$
Now since we can write $\partial\Sigma = c_{1} \sqcup \cdots \sqcup c_{n}$, we have an identification
$$ X(\partial\Sigma) = X(c_{1}) \times \dots \times X(c_{n}) \simeq \Abb^{n} $$
given by taking a local system on the disjoint union $\partial\Sigma$ of $n$ circles to its sequence of traces along the curves. The morphism $\textup{tr}_{\partial\Sigma}$ above may be viewed as an assignment to each $\rho \in X(\Sigma)$ its sequence of traces $\textup{tr}\,\rho(c_{1}),\dots,\textup{tr}\,\rho(c_{n})$. The fibers of $\textup{tr}_{\partial\Sigma}$ for $k \in \Abb^{n}$ will be denoted $X_{k} = X_{k}(\Sigma)$. Each $X_{k}$ is often called a \textit{relative character variety} in the literature. If $\Sigma$ is a surface of type $(g,n)$ satisfying $3g + n - 3 > 0$, then the relative character variety $X_{k}(\Sigma)$ is an irreducible algebraic variety of dimension $6g + 2n - 6$.

Given a fixed surface $\Sigma$, a subset $K \subseteq X(\partial\Sigma,\Cbb)$, and a subset $A \subseteq \Cbb$, we shall denote by
$$ X_{K}(A) = X_{K}(\Sigma,A) := X_{K}(\Sigma)(A) $$ 
the set of all $\rho \in X(\Sigma,\Cbb)$ such that $\textup{tr}_{\partial\Sigma}(\rho) \in K$ and $\textup{tr}_{a}(\rho) \in A$ for every essential curve $a \subset \Sigma$. By \cite[Lemma 2.5]{Wha20}, there is no risk of ambiguity with this notation, i.e., $X_{k}$ has a model over $A$ and $X_{k}(A)$ recovers the set of $A$-valued points of $X_{k}$ in the sense of algebraic geometry.

\subsection{Markoff-type cubic surfaces}
Now we give a description of the moduli spaces $X_{k}(\Sigma)$ for $(g,n) = (1,1)$ and $(0,4)$. These cases are special since each $X_{k}$ is an affine cubic algebraic surface with an explicit equation.

\begin{enumerate}
	\item[(1)] Let $\Sigma$ be a surface of type $(g,n) = (1,1)$, i.e. a one holed torus. Let $(\alpha,\beta,\gamma)$ be an \textit{optimal sequence of generators} for $\pi_{1}(\Sigma)$, as given in \cite[Definition 2.1]{Wha20}. By Example 2.2, we have an identification $X(\Sigma) \simeq \Abb^{3}$. From the trace relations in $\textup{SL}_{2}$, we obtain that 
	\[
	\begin{aligned}
	\textup{tr}_{\gamma} &= \textup{tr}_{\alpha\beta\alpha^{-1}\beta^{-1}} = \textup{tr}_{\alpha\beta\alpha^{-1}}\textup{tr}_{\beta^{-1}} - \textup{tr}_{\alpha\beta_\alpha^{-1}\beta}\\ &= \textup{tr}_{\beta}^{2} - \textup{tr}_{\alpha\beta}\textup{tr}_{\alpha^{-1}\beta} + \textup{tr}_{\alpha\alpha} = \textup{tr}_{\beta}^{2} - \textup{tr}_{\alpha\beta}(\textup{tr}_{\alpha^{-1}}\textup{tr}_{\beta} - \textup{tr}_{\alpha\beta}) + \textup{tr}_{\alpha}^{2} - \textup{tr}_{1}\\ &= \textup{tr}_{\alpha}^{2} + \textup{tr}_{\beta}^{2} + \textup{tr}_{\alpha\beta}^{2} - \textup{tr}_{\alpha}\textup{tr}_{\beta}\textup{tr}_{\alpha\beta} - 2.  
	\end{aligned}
	\]
	
	Writing $(x,y,z) = (\textup{tr}_{\alpha},\textup{tr}_{\beta},\textup{tr}_{\alpha\beta})$ so that each of the variables $x$, $y$, and $z$ corresponds to an essential curve on $\Sigma$ as depicted in \cite[Figure $2$]{Wha20}, the moduli space $X_{k} \subset X$ has an explicit presentation as an affine cubic algebraic surface in $\Abb^{3}_{x,y,z}$ with the equation
$$ x^{2} + y^{2} + z^{2} - xyz - 2 = k. $$
These are exactly the Markoff surfaces as studied in the series of papers \cite{GS22}, \cite{LM20}, and \cite{CTWX20} with $m = k+2$.	
	
	\item[(2)] Let $\Sigma$ be a surface of type $(g,n) = (0,4)$, i.e. a four holed sphere. Let $(\gamma_{1},\gamma_{2},\gamma_{3},\gamma_{4})$ be an optimal sequence of generators for $\pi_{1}(\Sigma)$. Set $$ (x,y,z) = (\textup{tr}_{\gamma_{1}\gamma_{2}}, \textup{tr}_{\gamma_{2}\gamma_{3}}, \textup{tr}_{\gamma_{1}\gamma_{3}}) $$
so that each of the variables corresponds to an essential curve on $\Sigma$. By Example 2.2, for $k = (k_{1},k_{2},k_{3},k_{4}) \in \Abb^{4}(\Cbb)$ the relative character variety $X_{k} = X_{k}(\Sigma)$ is an affine cubic algebraic surface in $\Abb^{3}_{x,y,z}$ given by the equation 
\begin{equation}
x^{2} + y^{2} + z^{2} + xyz = ax + by + cz + d,
\end{equation}
where 
\begin{equation}
	\begin{cases*}
	a = k_{1}k_{2} + k_{3}k_{4} \\
	b = k_{1}k_{4} + k_{2}k_{3} \\
	c = k_{1}k_{3} + k_{2}k_{4}
	\end{cases*}
	\hspace{0.5cm} \textup{and} \hspace{0.5cm} d = 4 - \sum_{i=1}^{4} k_{i}^{2} - \prod_{i=1}^{4} k_{i}. 
\end{equation}
These are the Markoff-type cubic surfaces that we are going to study in this paper.
\end{enumerate}

\section{The Brauer group of Markoff-type cubic surfaces}
Our main interest is in the \textbf{second} Markoff-type cubic surfaces defined by $(1)$. We are now going to give some explicit computations on the Brauer group of these surfaces. First of all, let us recall some basic definitions and results on the Brauer group of varieties over a field.

Let $k$ be an arbitrary field. Recall that for a variety $X$ over $k$ there is a natural filtration on the Brauer group
$$ \textup{Br}_{0}\,X \subset \textup{Br}_{1}\,X \subset \textup{Br}\,X $$ which is defined as follows.
\begin{definition}
Let
$$ \textup{Br}_{0}\,X = \textup{Im}[\textup{Br}\,k \rightarrow \textup{Br}\,X], \hspace{0.5cm} \textup{Br}_{1}\,X = \textup{Ker}[\textup{Br}\,X \rightarrow \textup{Br}\,\overline{X}]. $$
The subgroup $\textup{Br}_{1}\,X \subset \textup{Br}\,X$ is the \textbf{algebraic} Brauer group of $X$ and the quotient $\textup{Br}\,X/\textup{Br}_{1}\,X$ is the \textbf{transcendental} Brauer group of $X$.
\end{definition}
From the Hochschild--Serre spectral sequence, we have the following spectral sequence:
\begin{equation*}
E_{2}^{pq} = \textup{H}^{p}_{\textup{\'et}}(k, \textup{H}_{\textup{\'et}}^{q}(\overline{X},\Gmbb)) \Longrightarrow \textup{H}^{p+q}_{\textup{\'et}}(X,\Gmbb),
\end{equation*}
which is contravariantly functorial in the $k$-variety $X$. It gives rise to the functorial exact sequence of terms of low degree:
\begin{equation*}
\begin{aligned}
    0 &\longrightarrow \textup{H}^{1}(k,\overline{k}[X]^{\times}) \longrightarrow \textup{Pic}\,X \longrightarrow \textup{Pic}\,\overline{X}^{G_{k}} \longrightarrow \textup{H}^{2}(k, \overline{k}[X]^{\times}) \longrightarrow \textup{Br}_{1}\,X\\
    &\longrightarrow \textup{H}^{1}(k, \textup{Pic}\,\overline{X}) \longrightarrow \textup{Ker}[\textup{H}^{3}(k, \overline{k}[X]^{\times}) \rightarrow \textup{H}^{3}_{\textup{\'et}}(X,\Gmbb)]. 
\end{aligned}
\end{equation*}
Let $X$ be a variety over a field $k$ such that $\overline{k}[X]^{\times} = \overline{k}^{\times}$. By Hilbert’s Theorem 90 we have $\textup{H}^{1}(k, \overline{k}^{\times}) = 0$, then by the above sequence there is an exact sequence
\begin{equation*}
\begin{aligned}
    0 &\longrightarrow \textup{Pic}\,X \longrightarrow \textup{Pic}\,\overline{X}^{G_{k}} \longrightarrow \textup{Br}\,k \longrightarrow \textup{Br}_{1}\,X \\
	&\longrightarrow \textup{H}^{1}(k, \textup{Pic}\,\overline{X}) \longrightarrow \textup{Ker}[\textup{H}^{3}(k, \overline{k}^{\times}) \rightarrow \textup{H}^{3}_{\textup{\'et}}(X,\Gmbb)].
\end{aligned}
\end{equation*}
This sequence is also contravariantly functorial in $X$.

\begin{remark}
Let $X$ be a variety over a field $k$ such that $\overline{k}[X]^{\times} = \overline{k}^{\times}$. This assumption $\overline{k}[X]^{\times} = \overline{k}^{\times}$ holds for any proper, geometrically connected and geometrically reduced $k$-variety $X$.
\begin{enumerate} 
\item[(1)] If $X$ has a $k$-point, which defined a section of the structure morphism $X \rightarrow \textup{Spec}\,k$, then each of the maps $\textup{Br}\,k \longrightarrow \textup{Br}_{1}\,X$ and $\textup{H}^{3}(k, \overline{k}^{\times}) \rightarrow \textup{H}^{3}_{\textup{\'et}} (X, \Gmbb)$ has a retraction, hence is injective. (Then $\textup{Pic}\,X \longrightarrow \textup{Pic}\,\overline{X}^{G_{k}}$ is an isomorphism.) Therefore, we have an isomorphism
$$ \textup{Br}_{1}\,X/\textup{Br}\,k \cong \textup{H}^{1}(k, \textup{Pic}\,\overline{X}). $$ 

\item[(2)] If $k$ is a number field, then $\textup{H}^{3}(k, \overline{k}^{\times}) = 0$ (see \cite{CF67}, Chapter VII, Section 11.4, p. 199). Thus for a variety $X$ over a number field $k$ such that $\overline{k}[X]^{\times} = \overline{k}^{\times}$, we have an isomorphism 
$$ \textup{Br}_{1}\,X/\textup{Br}_{0}\,X \cong \textup{H}^{1}(k, \textup{Pic}\,\overline{X}). $$
\end{enumerate}
\end{remark}

\subsection{Geometry of affine cubic surfaces}
In this section, we study the geometry of affine cubic surfaces with special regards to the Brauer group. By an \textit{affine cubic surface}, we mean an affine surface of the form 
$$ U : f(u_{1}, u_{2}, u_{3}) = 0 $$
where $f$ is a polynomial of degree of 3. The closure of $U$ in $\Pbb^{3}$ is a cubic surface $X$. The complement $H = X\ \backslash\ U$ is a hyperplane section on $S$. Much of the geometry of $U$ can be understood in terms of the geometry of $X$ and $H$, especially in the case of Markoff-type cubic surfaces. There has been already much work on the Brauer groups of affine cubic surfaces when the hyperplane section $H$ is \textit{smooth}, for example see \cite{CTW12}. Here we shall be interested in the case where the hyperplane section $H$ is \textit{singular}; in particular, we focus on the case where $H$ is given by 3 coplanar lines. All results here are proven in either \cite{CTWX20} or \cite{LM20}.

We begin with an important result for cubic surfaces over an algebraically closed field.

\begin{proposition} \cite[Proposition 2.2]{CTWX20} 
	Let $X \subset \Pbb^{3}_{k}$ be a smooth projective cubic surface over a field $k$ of characteristic zero. Suppose a plane $\Pbb^{2}_{k} \subset \Pbb^{3}_{k}$ cuts out on $\overline{X}$ three distinct lines $L_{1}, L_{2}, L_{3}$ over $\overline{k}$. Let $U \subset X$ be the complement of this plane. Then the natural map $\overline{k}^{\times} \rightarrow \overline{k}[U]^{\times}$ is an isomorphism of Galois modules and the natural sequence
	$$ 0 \longrightarrow \bigoplus_{i=1}^{3} \Zbb L_{i} \longrightarrow \textup{Pic}\,\overline{X} \longrightarrow \textup{Pic}\,\overline{U} \longrightarrow 0 $$
	is an exact sequence of Galois lattices.
\end{proposition}

As $\text{Pic}\,\overline{U}$ is torsion free, we have the following result for the algebraic Brauer group, using the computation by $\textsc{Magma}$.

\begin{proposition} \cite[Proposition 2.5]{LM20}
Let $X$ be a smooth projective cubic surface over a field $k$ of characteristic 0. Let $H \subset S$ be a hyperplane section which is the union of 3 distinct lines $L_{1}, L_{2}, L_{3}$ and let $U = X\ \backslash\ H$. Then $\textup{Pic}\,\overline{U}$ is torsion free and $\textup{Br}_{1}\,U/\textup{Br}_{0}\,U \cong \textup{H}^{1}(k, \textup{Pic}\,\overline{U})$ is isomorphic to one of the following groups:
$$ 0, \Zbb/4\Zbb, \Zbb/2\Zbb \times \Zbb/4\Zbb, (\Zbb/2\Zbb)^{r} \hspace{0.5cm} (r = 1, 2, 3, 4). $$
\end{proposition}

For the transcendental Brauer group, from the discussion on page 140 of \cite{CTS21}, note that $\text{Br}_{1}\,U = \text{Ker}(\text{Br}\,U \rightarrow \text{Br}\,\overline{U}^{G_{k}})$ so we have $\text{Br}\,U/\text{Br}_{1}\,U \subset \text{Br}\,\overline{U}^{G_{k}}$.

\begin{proposition} \cite[Proposition 2.1]{CTWX20} \cite[Proposition 2.4]{LM20}
Let $X$ be a smooth projective cubic surface over a field $k$ of characteristic 0. Suppose that $U$ is an open subset of $X$ such that $X\ \backslash\ U$ is the union of three distinct $k$-lines, by which we mean a smooth projective curve isomorphic to $\Pbb^{1}_{k}$. Suppose any two lines intersect each another transversely in one point, and that the three intersection points are distinct. Let $L$ be one of the three lines and $V \subset L$ be the complement of the 2 intersection points of $L$ with the other two lines. Then the residue map
$$ \partial_{L} : \textup{Br}\,\overline{k}(X) \rightarrow \textup{H}^{1}(\overline{k}(L), \Qbb/\Zbb) $$
induces a $G_{k}$-isomorphism
$$ \textup{Br}\,\overline{U} \simeq \textup{H}^{1}(\overline{V}, \Qbb/\Zbb) \simeq \textup{H}^{1}(\overline{\mathbb{G}}_{m}, \Qbb/\Zbb) \simeq \Qbb/\Zbb(-1). $$
In particular, if $k$ contains no non-trivial roots of unity then
$$ \textup{Br}\,\overline{U}^{G_{k}} = \Zbb/2\Zbb. $$
\end{proposition}

\begin{lemma} \cite[Lemma 2.4]{CTWX20}
Let $k$ be a field of characteristic 0. Let $G_{k} = \textup{Gal}(\overline{k}/k)$. Then $\Qbb/\Zbb(-1)^{G_{k}}$ is (noncanonically) isomorphic to $\mu_{\infty}(k)$, the group of roots of unity in $k$. 
\end{lemma}

We end this section by the following result which applies to number fields and more generally to function fields of varieties over number fields.

\begin{corollary} \cite[Corollary 2.3]{CTWX20}
Let $k$ be a field of characteristic 0 such that in any finite field extension there are only finitely many roots of unity. Let $X \subset \Pbb^{3}_{k}$ be a smooth projective cubic surface over $k$. Suppose that a plane cuts out on $X$ three nonconcurrent lines. Let $U \subset X$ be the complement of the plane section. Then the quotient $\textup{Br}\,U/\textup{Br}_{0}\,U$ is finite.
\end{corollary}

\subsection{The geometric Picard group and algebraic Brauer group}
Using the equations, we can compute explicitly the algebraic Brauer group of the Markoff-type cubic surfaces in question. First, we have the following important result.

\begin{lemma}
Let $K$ be a number field and let $X \subset \Pbb^{3}_{K}$ be a cubic surface defined by the equation 
$$ t(x^{2}+y^{2}+z^{2}) + xyz = t^{2}(ax+by+cz) + dt^{3}, $$
where $a,b,c,d$ are defined by $(2)$ for some $k = (k_{1},k_{2},k_{3},k_{4}) \in \Abb^{4}(K)$. Then $X$ is singular if and only if we are in one of the following cases: 
\begin{enumerate}
	\item[$\bullet$] $\Delta(k) = 0$ where $k = (k_{1},k_{2},k_{3},k_{4}) \in \Abb^{4}(K)$ and 
	$$ \Delta(k) = (2(k_{1}^{2}+k_{2}^{2}+k_{3}^{2}+k_{4}^{2}) - k_{1}k_{2}k_{3}k_{4} - 16)^{2} - (4-k_{1}^{2})(4-k_{2}^{2})(4-k_{3}^{2})(4-k_{4}^{2}), $$	
	\item[$\bullet$] at least one of the parameters $k_{1},k_{2},k_{3},k_{4}$ equals $\pm 2$.
\end{enumerate} 
If $k$ satisfies none of those two conditions and $[E:K] = 16$ where 
$$ E := K(\sqrt{k_{1}^{2}-4},\sqrt{k_{2}^{2}-4},\sqrt{k_{3}^{2}-4},\sqrt{k_{4}^{2}-4}), $$
then the 27 lines on the smooth cubic surface $\overline{X}$ are defined over $E$ by the following equations 
$$ L_{1}: x=t=0; \hspace{1cm} L_{2}: y=t=0; \hspace{1cm} L_{3}: z=t=0 $$
and 
\begin{enumerate} 
	\item $\ell_{1}(\epsilon,\delta): x = \ds\frac{(k_{1}k_{2} + \epsilon\delta\sqrt{(k_{1}^{2}-4)(k_{2}^{2}-4)})}{2} t,\\ y = -\ds\frac{(k_{1}+\epsilon\sqrt{k_{1}^{2}-4})(k_{2}+\delta\sqrt{k_{2}^{2}-4})}{4} z - \ds\frac{c-b\frac{(k_{1}+\epsilon\sqrt{k_{1}^{2}-4})(k_{2}+\delta\sqrt{k_{2}^{2}-4})}{4}}{\frac{\delta k_{1}\sqrt{k_{2}^{2}-4}+\epsilon k_{2}\sqrt{k_{1}^{2}-4}}{2}} t;$	
	
	\item $\ell_{2}(\epsilon,\delta): y = \ds\frac{(k_{1}k_{4} + \epsilon\delta\sqrt{(k_{1}^{2}-4)(k_{4}^{2}-4)})}{2} t,\\ z = -\ds\frac{(k_{1}+\epsilon\sqrt{k_{1}^{2}-4})(k_{4}+\delta\sqrt{k_{4}^{2}-4})}{4} x - \ds\frac{a-c\frac{(k_{1}+\epsilon\sqrt{k_{1}^{2}-4})(k_{4}+\delta\sqrt{k_{4}^{2}-4})}{4}}{\frac{\delta k_{1}\sqrt{k_{4}^{2}-4}+\epsilon k_{4}\sqrt{k_{1}^{2}-4}}{2}} t;$
	
	\item $\ell_{3}(\epsilon,\delta): z = \ds\frac{(k_{1}k_{3} + \epsilon\delta\sqrt{(k_{1}^{2}-4)(k_{3}^{2}-4)})}{2} t,\\ y = -\ds\frac{(k_{1}+\epsilon\sqrt{k_{1}^{2}-4})(k_{3}+\delta\sqrt{k_{3}^{2}-4})}{4} x - \ds\frac{a-b\frac{(k_{1}+\epsilon\sqrt{k_{1}^{2}-4})(k_{3}+\delta\sqrt{k_{3}^{2}-4})}{4}}{\frac{\delta k_{1}\sqrt{k_{3}^{2}-4}+\epsilon k_{3}\sqrt{k_{1}^{2}-4}}{2}} t;$
	
	\item $\ell_{4}(\epsilon,\delta): x = \ds\frac{(k_{3}k_{4} + \epsilon\delta\sqrt{(k_{3}^{2}-4)(k_{4}^{2}-4)})}{2} t,\\ y = -\ds\frac{(k_{3}+\epsilon\sqrt{k_{3}^{2}-4})(k_{4}+\delta\sqrt{k_{4}^{2}-4})}{4} z - \ds\frac{c-b\frac{(k_{3}+\epsilon\sqrt{k_{3}^{2}-4})(k_{4}+\delta\sqrt{k_{4}^{2}-4})}{4}}{\frac{\delta k_{3}\sqrt{k_{4}^{2}-4}+\epsilon k_{4}\sqrt{k_{3}^{2}-4}}{2}} t;$
	
	\item $\ell_{5}(\epsilon,\delta): y = \ds\frac{(k_{2}k_{3} + \epsilon\delta\sqrt{(k_{2}^{2}-4)(k_{3}^{2}-4)})}{2} t,\\ z = -\ds\frac{(k_{2}+\epsilon\sqrt{k_{2}^{2}-4})(k_{3}+\delta\sqrt{k_{3}^{2}-4})}{4} x - \ds\frac{a-c\frac{(k_{2}+\epsilon\sqrt{k_{2}^{2}-4})(k_{3}+\delta\sqrt{k_{3}^{2}-4})}{4}}{\frac{\delta k_{2}\sqrt{k_{3}^{2}-4}+\epsilon k_{3}\sqrt{k_{2}^{2}-4}}{2}} t;$
	
	\item $\ell_{6}(\epsilon,\delta): z = \ds\frac{(k_{2}k_{4} + \epsilon\delta\sqrt{(k_{2}^{2}-4)(k_{4}^{2}-4)})}{2} t,\\ y = -\ds\frac{(k_{2}+\epsilon\sqrt{k_{2}^{2}-4})(k_{4}+\delta\sqrt{k_{4}^{2}-4})}{4} x - \ds\frac{a-b\frac{(k_{2}+\epsilon\sqrt{k_{2}^{2}-4})(k_{4}+\delta\sqrt{k_{4}^{2}-4})}{4}}{\frac{\delta k_{2}\sqrt{k_{4}^{2}-4}+\epsilon k_{4}\sqrt{k_{2}^{2}-4}}{2}} t$
	\end{enumerate}
with $\epsilon = \pm 1$ and $\delta = \pm 1$. Furthermore, we have the intersection numbers $$ \ell_{i}(\epsilon,\delta).\ell_{j}(\epsilon,\delta) = 0 $$ for any pair $(\epsilon,\delta)$, for all $1 \leq i \not= j \leq 6$.
\end{lemma}

\begin{proof}
The necessary and sufficient condition for the affine open surface $U = X \setminus \{t = 0\}$ to be singular is proven in \cite[Theorem 3.7]{CL09}. It is easy to verify that there is no singular point at infinity on the projective surface $X$.

Now without loss of generality, we consider the system of equations
\begin{equation*}
\begin{cases*}
y = \alpha_{1}x + \alpha_{2}t \\
z = \beta_{1}x + \beta_{2}t
\end{cases*}
\end{equation*}
and put them in the original equation of the cubic surfaces to solve $\alpha_{i},\beta_{i}$ for $i=1,2$. We can work similarly for $(z,x)$ and $(x,y)$ to find all the given equations of the 27 lines.
\end{proof}

Now given the data of the lines, we can compute directly the algebraic Brauer group of the Markoff-type cubic surfaces in question.

\begin{proposition}
Let $K$ be a number field. Let $X \subset \Pbb^{3}_{K}$ be a cubic surface defined by the equation 
\begin{equation}
t(x^{2}+y^{2}+z^{2}) + xyz = t^{2}(ax+by+cz) + dt^{3},
\end{equation}
where $a,b,c,d$ are defined by $(2)$ for some $k = (k_{1},k_{2},k_{3},k_{4}) \in \Abb^{4}(K)$. Assume that $X$ is smooth over $K$ and $[E : K] = 16$, then 
$$ \textup{Br}\,X/\textup{Br}_{0}\,X = \textup{Br}_{1}\,X/\textup{Br}_{0}\,X \cong \Zbb/2\Zbb. $$
\end{proposition}

\begin{proof}
Since $X$ is geometrically rational, one has $\textup{Br}\,X = \textup{Br}_{1}\,X$. By taking $x=t=0$ for instance, one clearly has $X(K) \not= \emptyset$, so $\textup{Br}_{0}\,X = \textup{Br}\,K$. Since $K$ is a number field, by the Hochschild--Serre spectral sequence, we have an isomorphism 
$$ \textup{Br}_{1}\,X/\textup{Br}_{0}\,X \simeq \textup{H}^{1}(K,\textup{Pic}\,\overline{X}). $$

By the above lemma, we can easily verify that the six lines $\ell_{1}(1,1)$, $\ell_{1}(1,-1)$, $\ell_{3}(-1,1)$, $\ell_{4}(-1,-1)$, $\ell_{4}(-1,1)$, and $L_{2}$ on the cubic surface $\overline{X}$ are skew to each other, hence they may be simultaneously blown down to $\Pbb^{2}$ by \cite{Har77}, Chapter V, Proposition 4.10. For the sake of simplicity, here we shall write these six lines respectively as $\ell_{i}$ for $1 \leq i \leq 6$. The class $\omega$ of the canonical divisor on $\overline{X}$ is equal to $-3\ell + \Sigma_{i=1}^{6} \ell_{i}$, where $\ell$ is the inverse image of the class of lines in $\Pbb^{2}$. By \cite[Chapter V, Proposition 4.8]{Har77}, the classes $\ell, \ell_{i}, i = \overline{1,6}$ form a basis of $\textup{Pic}\,\overline{X}$, and we have the following intersection properties: $(\ell.\ell) = 1, (\ell.\ell_{i}) = 0$ for $1 \leq i \leq 6$. 

Since $(L_{1}.\ell_{3}) = 0, (L_{1}.\ell_{i}) = 1, i \not= 3$; $(L_{3}.\ell_{3}) = (L_{3}.\ell_{6}) = 1, (L_{3}.\ell_{i}) = 0, i \not= 3,6$; and $L_{2} = \ell_{6}$, one concludes that 
\begin{equation}
L_{1} = 2\ell - \Sigma_{i \not= 3} \ell_{i}, \hspace{1cm} L_{3} = \ell - \ell_{3} - \ell_{6}
\end{equation} 
in $\textup{Pic}\,\overline{X}$ by \cite[Chapter V, Proposition 4.8 (e)]{Har77}.

Now we consider the action of the Galois group $G := \textup{Gal}(E/K)$ on $\textup{Pic}\,\overline{X}$. One clearly has $G \cong \langle \sigma_{1} \rangle \times \langle \sigma_{2} \rangle \times \langle \sigma_{3} \rangle \times \langle \sigma_{4} \rangle$, where 
$$ \sigma_{i}(\sqrt{k_{i}^{2}-4}) = -\sqrt{k_{i}^{2}-4} \hspace{0.5cm} \textup{and} \hspace{0.5cm} \sigma_{i}(\sqrt{k_{j}^{2}-4}) = \sqrt{k_{j}^{2}-4}, 1 \leq i \not= j \leq 4. $$
We have the following intersection numbers, noting that $\sigma_{2}(\ell_{1}) = \ell_{2}$, $\sigma_{2}(\ell_{2}) = \ell_{1}$, $\sigma_{4}(\ell_{4}) = \ell_{5}$, and $\sigma_{4}(\ell_{5}) = \ell_{4}$:

\begin{equation}
\begin{cases}
(\sigma_{1}(\ell_{1}).\ell_{1}) = (\ell_{1}(-1,1).\ell_{1}(1,1)) = 0\\
(\sigma_{1}(\ell_{1}).\ell_{2}) = (\ell_{1}(-1,1).\ell_{1}(1,-1)) = 1\\
(\sigma_{1}(\ell_{1}).\ell_{3}) = (\ell_{1}(-1,1).\ell_{3}(-1,1)) = 1\\
(\sigma_{1}(\ell_{1}).\ell_{4}) = (\ell_{1}(-1,1).\ell_{4}(-1,-1)) = 0\\
(\sigma_{1}(\ell_{1}).\ell_{5}) = (\ell_{1}(-1,1).\ell_{4}(-1,1)) = 0\\
(\sigma_{1}(\ell_{1}).\ell_{6}) = (\ell_{1}(-1,1).L_{2}) = 0,
\end{cases}
\end{equation}

\begin{equation}
\begin{cases}
(\sigma_{1}(\ell_{2}).\ell_{1}) = (\ell_{1}(-1,-1).\ell_{1}(1,1)) = 1\\
(\sigma_{1}(\ell_{2}).\ell_{2}) = (\ell_{1}(-1,-1).\ell_{1}(1,-1)) = 0\\
(\sigma_{1}(\ell_{2}).\ell_{3}) = (\ell_{1}(-1,-1).\ell_{3}(-1,1)) = 1\\
(\sigma_{1}(\ell_{2}).\ell_{4}) = (\ell_{1}(-1,-1).\ell_{4}(-1,-1)) = 0\\
(\sigma_{1}(\ell_{2}).\ell_{5}) = (\ell_{1}(-1,-1).\ell_{4}(-1,1)) = 0\\
(\sigma_{1}(\ell_{2}).\ell_{6}) = (\ell_{1}(-1,-1).L_{2}) = 0,
\end{cases}
\end{equation}

\begin{equation}
\begin{cases}
(\sigma_{1}(\ell_{3}).\ell_{1}) = (\ell_{3}(1,1).\ell_{1}(1,1)) = 1\\
(\sigma_{1}(\ell_{3}).\ell_{2}) = (\ell_{3}(1,1).\ell_{1}(1,-1)) = 1\\
(\sigma_{1}(\ell_{3}).\ell_{3}) = (\ell_{3}(1,1).\ell_{3}(-1,1)) = 0\\
(\sigma_{1}(\ell_{3}).\ell_{4}) = (\ell_{3}(1,1).\ell_{4}(-1,-1)) = 0\\
(\sigma_{1}(\ell_{3}).\ell_{5}) = (\ell_{3}(1,1).\ell_{4}(-1,1)) = 0\\
(\sigma_{1}(\ell_{3}).\ell_{6}) = (\ell_{3}(1,1).L_{2}) = 0;
\end{cases}
\end{equation}
and
\begin{equation}
\begin{cases}
(\sigma_{3}(\ell_{3}).\ell_{1}) = (\ell_{3}(-1,-1).\ell_{1}(1,1)) = 0\\
(\sigma_{3}(\ell_{3}).\ell_{2}) = (\ell_{3}(-1,-1).\ell_{1}(1,-1)) = 0\\
(\sigma_{3}(\ell_{3}).\ell_{3}) = (\ell_{3}(-1,-1).\ell_{3}(-1,1)) = 0\\
(\sigma_{3}(\ell_{3}).\ell_{4}) = (\ell_{3}(-1,-1).\ell_{4}(-1,-1)) = 1\\
(\sigma_{3}(\ell_{3}).\ell_{5}) = (\ell_{3}(-1,-1).\ell_{4}(-1,1)) = 1\\
(\sigma_{3}(\ell_{3}).\ell_{6}) = (\ell_{3}(-1,-1).L_{2}) = 0,
\end{cases}
\end{equation} 

\begin{equation}
\begin{cases}
(\sigma_{3}(\ell_{4}).\ell_{1}) = (\ell_{4}(1,-1).\ell_{1}(1,1)) = 0\\
(\sigma_{3}(\ell_{4}).\ell_{2}) = (\ell_{4}(1,-1).\ell_{1}(1,-1)) = 0\\
(\sigma_{3}(\ell_{4}).\ell_{3}) = (\ell_{4}(1,-1).\ell_{3}(-1,1)) = 1\\
(\sigma_{3}(\ell_{4}).\ell_{4}) = (\ell_{4}(1,-1).\ell_{4}(-1,-1)) = 0\\
(\sigma_{3}(\ell_{4}).\ell_{5}) = (\ell_{4}(1,-1).\ell_{4}(-1,1)) = 1\\
(\sigma_{3}(\ell_{4}).\ell_{6}) = (\ell_{4}(1,-1).L_{2}) = 0,
\end{cases}
\end{equation} 

\begin{equation}
\begin{cases}
(\sigma_{3}(\ell_{5}).\ell_{1}) = (\ell_{4}(1,1).\ell_{1}(1,1)) = 0\\
(\sigma_{3}(\ell_{5}).\ell_{2}) = (\ell_{4}(1,1).\ell_{1}(1,-1)) = 0\\
(\sigma_{3}(\ell_{5}).\ell_{3}) = (\ell_{4}(1,1).\ell_{3}(-1,1)) = 1\\
(\sigma_{3}(\ell_{5}).\ell_{4}) = (\ell_{4}(1,1).\ell_{4}(-1,-1)) = 1\\
(\sigma_{3}(\ell_{5}).\ell_{5}) = (\ell_{4}(1,1).\ell_{4}(-1,1)) = 0\\
(\sigma_{3}(\ell_{5}).\ell_{6}) = (\ell_{4}(1,1).L_{2}) = 0.
\end{cases}
\end{equation} 

Hence, we obtain
\begin{equation}
\begin{cases}
\sigma_{1}(\ell_{1}) = \ell - \ell_{2} - \ell_{3}\\
\sigma_{1}(\ell_{2}) = \ell - \ell_{1} - \ell_{3}\\
\sigma_{1}(\ell_{3}) = \ell - \ell_{1} - \ell_{2}\\
\sigma_{3}(\ell_{3}) = \ell - \ell_{4} - \ell_{5}\\
\sigma_{3}(\ell_{4}) = \ell - \ell_{3} - \ell_{5}\\
\sigma_{3}(\ell_{5}) = \ell - \ell_{3} - \ell_{4}
\end{cases}
\end{equation}
in $\textup{Pic}\,\overline{X}$ by \cite[Chapter V, Proposition 4.9]{Har77}. As a result, we have 
\begin{equation}
\begin{cases}
\sigma_{1}(\ell) = 2\ell - \ell_{1} - \ell_{2} - \ell_{3}\\
\sigma_{3}(\ell) = 2\ell - \ell_{3} - \ell_{4} - \ell_{5},
\end{cases}
\end{equation}
and clearly $\sigma_{2}(\ell) = \sigma_{4}(\ell) = \ell$. Then 
\begin{equation}
\begin{cases}
\textup{Ker}(1+\sigma_{1}) = \langle \ell-\ell_{1}-\ell_{2}-\ell_{3} \rangle\\
\textup{Ker}(1+\sigma_{2}) = \langle \ell_{1}-\ell_{2} \rangle\\
\textup{Ker}(1+\sigma_{3}) = \langle \ell-\ell_{3}-\ell_{4}-\ell_{5} \rangle\\
\textup{Ker}(1+\sigma_{4}) = \langle \ell_{4}-\ell_{5} \rangle,\\
\end{cases}
\end{equation}

\begin{equation}
\begin{cases}
\textup{Ker}(1-\sigma_{1}) = \langle \ell-\ell_{1}, \ell-\ell_{2}, \ell-\ell_{3}, \ell_{4}, \ell_{5}, \ell_{6} \rangle\\
\textup{Ker}(1-\sigma_{2}) = \langle \ell, \ell_{1}+\ell_{2},\ell_{3},\ell_{4},\ell_{5},\ell_{6} \rangle\\
\textup{Ker}(1-\sigma_{3}) = \langle \ell-\ell_{3}, \ell-\ell_{4}, \ell-\ell_{5}, \ell_{1}, \ell_{2}, \ell_{6} \rangle\\
\textup{Ker}(1-\sigma_{4}) = \langle \ell, \ell_{1}, \ell_{2}, \ell_{3}, \ell_{4}+\ell_{5}, \ell_{6} \rangle,\\
\end{cases}
\end{equation}
and 
\begin{equation}
\begin{cases}
(1-\sigma_{1})\textup{Pic}\,\overline{X} = \langle \ell-\ell_{1}-\ell_{2}-\ell_{3} \rangle\\
(1-\sigma_{2})\textup{Pic}\,\overline{X} = \langle \ell_{1}-\ell_{2} \rangle\\
(1-\sigma_{3})\textup{Pic}\,\overline{X} = \langle \ell-\ell_{3}-\ell_{4}-\ell_{5} \rangle\\
(1-\sigma_{4})\textup{Pic}\,\overline{X} = \langle \ell_{4}-\ell_{5} \rangle.\\
\end{cases}
\end{equation}

Given a finite cyclic group $G = \langle \sigma \rangle$ and a $G$-module $M$, by \cite[Proposition 1.7.1]{NSW15}, recall that we have isomorphisms $\textup{H}^{1}(G,M) \cong \hat{\textup{H}}^{-1}(G,M)$, where the latter group is the quotient of $\prescript{}{N_{G}}{M}$, the set of elements of $M$ of norm $0$, by its subgroup $(1 - \sigma)M$.

By \cite[Proposition 1.6.7]{NSW15}, we have 
$$ \textup{H}^{1}(K,\textup{Pic}\,\overline{X}) = \textup{H}^{1}(G,\textup{Pic}\,\overline{X}), $$ where $G = \langle \sigma_{i}, 1 \leq i \leq 4 \rangle$. Then one has the following (inflation-restriction) exact sequence
$$ 0 \rightarrow \textup{H}^{1}(\langle \sigma_{1},\sigma_{2},\sigma_{3} \rangle, \textup{Pic}\,\overline{X}^{\langle \sigma_{4} \rangle}) \rightarrow \textup{H}^{1}(G, \textup{Pic}\,\overline{X}) \rightarrow \textup{H}^{1}(\langle \sigma_{4} \rangle, \textup{Pic}\,\overline{X}) = 0, $$
hence $\textup{H}^{1}(G, \textup{Pic}\,\overline{X}) \cong \textup{H}^{1}(\langle \sigma_{1},\sigma_{2},\sigma_{3} \rangle, \textup{Pic}\,\overline{X}^{\langle \sigma_{4} \rangle})$. Continuing as above, we have 
$$ 0 \rightarrow \textup{H}^{1}(\langle \sigma_{1},\sigma_{3} \rangle, \textup{Pic}\,\overline{X}^{\langle \sigma_{2}, \sigma_{4} \rangle}) \rightarrow \textup{H}^{1}(\langle \sigma_{1},\sigma_{2},\sigma_{3} \rangle, \textup{Pic}\,\overline{X}^{\langle \sigma_{4} \rangle}) \rightarrow \textup{H}^{1}(\langle \sigma_{2} \rangle, \textup{Pic}\,\overline{X}^{\langle \sigma_{4} \rangle}) = 0, $$
hence $\textup{H}^{1}(\langle \sigma_{1},\sigma_{2},\sigma_{3} \rangle, \textup{Pic}\,\overline{X}^{\langle \sigma_{4} \rangle}) \cong \textup{H}^{1}(\langle \sigma_{1},\sigma_{3} \rangle, \textup{Pic}\,\overline{X}^{\langle \sigma_{2}, \sigma_{4} \rangle})$. Now we are left with 
$$ 0 \rightarrow \textup{H}^{1}(\langle \sigma_{1} \rangle, \textup{Pic}\,\overline{X}^{\langle \sigma_{2}, \sigma_{3}, \sigma_{4} \rangle}) \rightarrow \textup{H}^{1}(\langle \sigma_{1},\sigma_{3} \rangle, \textup{Pic}\,\overline{X}^{\langle \sigma_{2}, \sigma_{4} \rangle}) \rightarrow \textup{H}^{1}(\langle \sigma_{3} \rangle, \textup{Pic}\,\overline{X}^{\langle \sigma_{2}, \sigma_{4} \rangle}) = 0, $$
hence $\textup{H}^{1}(\langle \sigma_{1},\sigma_{3} \rangle, \textup{Pic}\,\overline{X}^{\langle \sigma_{2}, \sigma_{4} \rangle}) \cong \textup{H}^{1}(\langle \sigma_{1} \rangle, \textup{Pic}\,\overline{X}^{\langle \sigma_{2}, \sigma_{3}, \sigma_{4} \rangle}) = \Zbb/2\Zbb$. Indeed, the last group can be computed as follows. We have 
$$ \textup{Pic}\,\overline{X}^{\langle \sigma_{2}, \sigma_{3}, \sigma_{4} \rangle} = \langle \ell_{1}+\ell_{2}, \ell-\ell_{3}, 2\ell-\ell_{4}-\ell_{5}, \ell_{6} \rangle. $$
Considering the action of $\sigma_{1}$ on this invariant group, we have
$$ \prescript{}{N_{\sigma_{1}}}{\textup{Pic}\,\overline{X}^{\langle \sigma_{2}, \sigma_{3}, \sigma_{4} \rangle}} = \textup{Ker}(1+\sigma_{1}) \cap \textup{Pic}\,\overline{X}^{\langle \sigma_{2}, \sigma_{3}, \sigma_{4} \rangle} = \langle \ell - \ell_{1} - \ell_{2} - \ell_{3} \rangle. $$
On the other hand, 
$$ (1-\sigma_{1})\textup{Pic}\,\overline{X}^{\langle \sigma_{2}, \sigma_{3}, \sigma_{4} \rangle} = [(1-\sigma_{1})\textup{Pic}\,\overline{X}] \cap \textup{Pic}\,\overline{X}^{\langle \sigma_{2}, \sigma_{3}, \sigma_{4} \rangle} = \langle 2(\ell - \ell_{1} - \ell_{2} - \ell_{3}) \rangle. $$
Given these results, we conclude that 
$$ \textup{H}^{1}(K,\textup{Pic}\,\overline{X}) = \textup{H}^{1}(G,\textup{Pic}\,\overline{X}) \cong \Zbb/2\Zbb. $$
\end{proof}

\begin{theorem}
Let $K$ be a number field. With the same notations as before, let $k \in \Abb^{4}(K)$ such that $[E : K] = 16$ where $E = K(\sqrt{k_{i}^{2}-4}, 1 \leq i \leq 4)$ and $X$ is smooth over $K$. Let $U$ be the affine cubic surface defined by the equation
$$ x^{2} + y^{2} + z^{2} + xyz = ax + by + cz + d, $$
where $a,b,c,d$ are defined by $(2)$ for some $k = (k_{1},k_{2},k_{3},k_{4}) \in \Abb^{4}(K)$. Then we have 
$$ \textup{Br}_{1}\,U/\textup{Br}_{0}\,U \cong \Zbb/2\Zbb $$
with a generator 
\[
\begin{aligned}
\mathcal{A} &= \textup{Cor}^{F_{1}}_{K} \left( x - \frac{k_{1}k_{2} + \sqrt{(k_{1}^{2}-4)(k_{2}^{2}-4)}}{2}, (k_{1}\sqrt{k_{2}^{2}-4} + k_{2}\sqrt{k_{1}^{2}-4})^{2} \right)\\ &= \textup{Cor}^{F_{3}}_{K} \left( y - \frac{k_{1}k_{4} + \sqrt{(k_{1}^{2}-4)(k_{4}^{2}-4)}}{2}, (k_{1}\sqrt{k_{4}^{2}-4} + k_{4}\sqrt{k_{1}^{2}-4})^{2} \right)\\ &= \textup{Cor}^{F_{2}}_{K} \left( z - \frac{k_{1}k_{3} + \sqrt{(k_{1}^{2}-4)(k_{3}^{2}-4)}}{2}, (k_{1}\sqrt{k_{3}^{2}-4} + k_{3}\sqrt{k_{1}^{2}-4})^{2} \right),
\end{aligned}
\]
where $F_{i} = K(\sqrt{(k_{1}^{2}-4)(k_{i+1}^{2}-4)})$ for $1 \leq i \leq 3$. Furthermore, we also have 
$$ \textup{Br}\,X/\textup{Br}_{0}\,X = \textup{Br}_{1}\,X/\textup{Br}_{0}\,X \cong \textup{Br}_{1}\,U/\textup{Br}_{0}\,U $$
with a generator 
\[
\begin{aligned}
\mathcal{A}_{0} &= \textup{Cor}^{F_{1}}_{K} \left( \frac{x}{t} - \frac{k_{1}k_{2} + \sqrt{(k_{1}^{2}-4)(k_{2}^{2}-4)}}{2}, (k_{1}\sqrt{k_{2}^{2}-4} + k_{2}\sqrt{k_{1}^{2}-4})^{2} \right)\\ &= \textup{Cor}^{F_{3}}_{K} \left( \frac{y}{t} - \frac{k_{1}k_{4} + \sqrt{(k_{1}^{2}-4)(k_{4}^{2}-4)}}{2}, (k_{1}\sqrt{k_{4}^{2}-4} + k_{4}\sqrt{k_{1}^{2}-4})^{2} \right)\\ &= \textup{Cor}^{F_{2}}_{K} \left( \frac{z}{t} - \frac{k_{1}k_{3} + \sqrt{(k_{1}^{2}-4)(k_{3}^{2}-4)}}{2}, (k_{1}\sqrt{k_{3}^{2}-4} + k_{3}\sqrt{k_{1}^{2}-4})^{2} \right)
\end{aligned}
\]
over $t \not= 0$.
\end{theorem}

\begin{proof}
We keep the notation as in the previous proposition. Then $\textup{Pic}\,\overline{U}$ is given by the following quotient group
$$ \textup{Pic}\,\overline{U} \cong \textup{Pic}\,\overline{X}/(\oplus_{i=1}^{3} \Zbb L_{i}) \cong (\oplus_{i=1}^{6} \Zbb \ell_{i} \oplus \Zbb \ell)/\langle 2\ell-\Sigma_{i \not= 3} \ell_{i}, \ell-\ell_{3}-\ell_{6}, \ell_{6}) \rangle \cong \oplus_{i=1}^{4} \Zbb [\ell_{i}] $$
by Proposition 3.1 and formula $(4)$. Here for any divisor $D \in \textup{Pic}\,\overline{X}$, denote by $[D]$ its image in $\textup{Pic}\,\overline{U}$. By Proposition 3.1, we also have $\overline{K}^{\times} = \overline{K}[U]^{\times}$. By the Hochschild--Serre spectral sequence, we have the following injective homomorphism 
$$ \textup{Br}_{1}\,U/\textup{Br}_{0}\,U \hookrightarrow \textup{H}^{1}(K,\textup{Pic}\,\overline{U}), $$
and in fact it is an \textbf{isomorphism} because over a number field $K$, we have $\textup{H}^{3}(K,\Gmbb) = 0$ from class field theory. Furthermore, the smooth compactification $X$ of $U$ has rational points, hence so does $U$, which comes from the fact that any smooth cubic surface over an infinite field $k$ is unirational over $k$ as soon as it has a $k$-point (see \cite{Kol02}), so we also have $\textup{Br}_{0}\,U = \textup{Br}\,K$.

Since $\textup{Pic}\,\overline{U}$ is free and $\textup{Gal}(\overline{K}/E)$ acts on $\textup{Pic}\,\overline{U}$ trivially, we obtain that $\textup{H}^{1}(K,\textup{Pic}\,\overline{U}) \cong \textup{H}^{1}(G,\textup{Pic}\,\overline{U})$ by \cite[Proposition 1.6.7]{NSW15}. Now in $\textup{Pic}\,\overline{U}$, as $[\ell_{6}] = 0$, $[\ell] = [\ell_{3}]$ and $2[\ell_{3}] = [\ell_{1}] + [\ell_{2}] + [\ell_{4}] + [\ell_{5}]$, we have the following equalities
\begin{equation}
\begin{cases}
\sigma_{1}([\ell_{1}]) = -[\ell_{2}], \sigma_{1}([\ell_{2}]) = -[\ell_{1}], \sigma_{1}([\ell_{3}]) = [\ell_{3}] - [\ell_{1}] - [\ell_{2}], \sigma_{1}([\ell_{4}]) = [\ell_{4}];\\
\sigma_{2}([\ell_{1}]) = [\ell_{2}], \sigma_{2}([\ell_{2}]) = [\ell_{1}], \sigma_{2}([\ell_{3}]) = [\ell_{3}], \sigma_{2}([\ell_{4}]) = [\ell_{4}];\\
\sigma_{3}([\ell_{1}]) = [\ell_{1}], \sigma_{3}([\ell_{2}]) = [\ell_{2}], \sigma_{3}([\ell_{3}]) = [\ell_{1}] + [\ell_{2}] - [\ell_{3}], \sigma_{3}([\ell_{4}]) = [\ell_{1}]+[\ell_{2}]+[\ell_{4}] - 2[\ell_{3}];\\
\sigma_{4}([\ell_{1}]) = [\ell_{1}], \sigma_{4}([\ell_{2}]) = [\ell_{2}], \sigma_{4}([\ell_{3}]) = [\ell_{3}], \sigma_{4}([\ell_{4}]) = 2[\ell_{3}] - [\ell_{1}] - [\ell_{2}] - [\ell_{4}].
\end{cases}
\end{equation}

Using the inflation-restriction sequences similarly as in the previous proposition, with $\textup{Pic}\,\overline{X}$ replaced by $\textup{Pic}\,\overline{U}$, we can compute that 
$$ \textup{H}^{1}(G,\textup{Pic}\,\overline{U}) \cong \textup{H}^{1}(\langle \sigma_{1} \rangle, \textup{Pic}\,\overline{U}^{\langle \sigma_{2},\sigma_{3},\sigma_{4} \rangle}) = \textup{H}^{1}(\langle \sigma_{1} \rangle, \langle [\ell_{1}] + [\ell_{2}] \rangle) = \ds \frac{\langle [\ell_{1}] + [\ell_{2}] \rangle}{\langle 2([\ell_{1}] + [\ell_{2}]) \rangle} \cong \Zbb/2\Zbb. $$

Now we produce concrete generators in $\textup{Br}_{1}\,U$ for $\textup{Br}_{1}\,U/\textup{Br}_{0}\,U$. Then $U$ is defined by the equation 
$$ x^{2} + y^{2} + z^{2} + xyz = ax + by + cz + d. $$
We will show that the following quaternion algebras in $\textup{Br}\,K(U)$ are non-constant elements of $\textup{Br}_{1}\,U$, and hence they are equal in $\textup{Br}_{1}\,U/\textup{Br}_{0}\,U$.
\begin{equation}
	\begin{cases}
	\mathcal{A}_{1} = \textup{Cor}^{F_{1}}_{K} \left( x - \frac{k_{1}k_{2} + \sqrt{(k_{1}^{2}-4)(k_{2}^{2}-4)}}{2}, (k_{1}\sqrt{k_{2}^{2}-4} + k_{2}\sqrt{k_{1}^{2}-4})^{2} \right), \\
	\mathcal{A}_{2} = \textup{Cor}^{F_{3}}_{K} \left( y - \frac{k_{1}k_{4} + \sqrt{(k_{1}^{2}-4)(k_{4}^{2}-4)}}{2}, (k_{1}\sqrt{k_{4}^{2}-4} + k_{4}\sqrt{k_{1}^{2}-4})^{2} \right), \\
	\mathcal{A}_{3} = \textup{Cor}^{F_{2}}_{K} \left( z - \frac{k_{1}k_{3} + \sqrt{(k_{1}^{2}-4)(k_{3}^{2}-4)}}{2}, (k_{1}\sqrt{k_{3}^{2}-4} + k_{3}\sqrt{k_{1}^{2}-4})^{2} \right).
	\end{cases} 
\end{equation}  
Indeed, it suffices to prove the claim for $\mathcal{A}_{1}$ and we only need to show that $\mathcal{A}_{1} \in \textup{Br}\,U$, since its formula implies that $\mathcal{A}_{1}$ becomes zero under the field extension $K \subset K(\sqrt{k_{1}^{2}-4}, \sqrt{k_{2}^{2}-4})$, i.e., it is algebraic. By Grothendieck's purity theorem (\cite[Theorem 6.8.3]{Poo17}), for any smooth integral variety $Y$ over a field $L$ of characteristic $0$, we have the exact sequence
$$ 0 \rightarrow \textup{Br}\,Y \rightarrow \textup{Br}\,L(Y) \rightarrow \oplus_{D \in Y^{(1)}} \textup{H}^{1}(L(D),\Qbb/\Zbb), $$
where the last map is given by the residue along the codimension-one point $D$. Therefore, to prove that our algebras come from a class in $\textup{Br}\,U$, it suffices to show that all their residues are trivial. We will show that
$$ \mathcal{A}_{1}' = \left( x - \frac{k_{1}k_{2} + \sqrt{(k_{1}^{2}-4)(k_{2}^{2}-4)}}{2}, (k_{1}\sqrt{k_{2}^{2}-4} + k_{2}\sqrt{k_{1}^{2}-4})^{2} \right) \in \textup{Br}\,U_{F_{1}} $$
so that its corestriction is a well-defined element over $K$. From the data of the 27 lines in Lemma 3.6 and the formula of $\mathcal{A}_{1}'$, any non-trivial residue of $\mathcal{A}_{1}'$ must occur along an irreducible component of the following divisor(s)
\[ 
\begin{aligned}
D_{1} : x &= \ds\frac{(k_{1}k_{2} + \sqrt{(k_{1}^{2}-4)(k_{2}^{2}-4)})}{2},\\ y &= -\ds\frac{(k_{1} \pm \sqrt{k_{1}^{2}-4})(k_{2} \pm \sqrt{k_{2}^{2}-4})}{4} z - \ds\frac{c-b\frac{(k_{1} \pm \sqrt{k_{1}^{2}-4})(k_{2} \pm \sqrt{k_{2}^{2}-4})}{4}}{\pm \frac{k_{1}\sqrt{k_{2}^{2}-4}+k_{2}\sqrt{k_{1}^{2}-4}}{2}}. 
\end{aligned}
\]
However, clearly in the function field of any such irreducible component, $(k_{1}\sqrt{k_{2}^{2}-4} + k_{2}\sqrt{k_{1}^{2}-4})^{2}$ is a square; standard formulae for residues in terms of the tame symbol \cite[Example 7.1.5, Proposition 7.5.1]{GS17} therefore show that $\mathcal{A}_{1}'$ is unramified, and hence $\mathcal{A}_{1} \in \textup{Br}\,U$. The residues of $\mathcal{A}_{1}$ at the lines $L_{1}, L_{2}, L_{3}$ which form the complement of $U$ in $X$ are easily seen to be trivial. One thus also has $\mathcal{A}_{1} \in \textup{Br}\,X$. This element is non-constant by the Faddeev exact sequence (\cite[Theorem 1.5.2]{CTS21}), since the residue of $\mathcal{A}_{1}'$ regarded as an element of $\textup{Br}\,F_{1}(x) = \textup{Br}\,F_{1}(\Pbb^{1})$ is nontrivial at the closed point $\left( x - \frac{k_{1}k_{2} + \sqrt{(k_{1}^{2}-4)(k_{2}^{2}-4)}}{2} \right)$ of $\Pbb^{1}_{F_{1}}$. Alternatively, by \cite[Corollary 11.3.5]{CTS21} on the Brauer group of conic bundles, the element $\mathcal{A}_{1}$ is indeed non-constant. Furthermore, it will contribute to the Brauer--Manin obstruction to strong approximation in the next section.

Finally, the fact that $\mathcal{A}_{1} = \mathcal{A}_{2} = \mathcal{A}_{3} = \mathcal{A}$ viewed as an element of $\textup{Br}_{1}\,U/\textup{Br}_{0}\,U \cong \Zbb/2\Zbb$ (by abuse of notation) and that $U$ is the open subset of $X$ defined by $t \not= 0$ give the desired generator $\mathcal{A}_{0}$ of $\textup{Br}\,X/\textup{Br}_{0}\,X \cong \Zbb/2\Zbb$.
\end{proof}

\subsection{The transcendental Brauer group}
We begin with a specific assumption. 

\begin{assumption}
Let $k = (k_{1},k_{2},k_{3},k_{4}) \in \Abb^{4}(\Zbb)$ and $E := \Qbb(\sqrt{k_{i}^{2}-4}, 1 \leq i \leq 4)$ such that $[E:\Qbb]=16$. For all $1 \leq i \not= j \leq 4$, assume that $|k_{i}| \geq 3$ such that $(k_{i} + \sqrt{k_{i}^{2}-4})(k_{j} + \sqrt{k_{j}^{2}-4})$ is not a square in $E$.
\end{assumption}

Now we compute the transcendental Brauer group in our particular case.
\begin{proposition}
Let $k = (k_{1},k_{2},k_{3},k_{4}) \in \Abb^{4}(\Zbb)$ satisfy Assumption 3.3. Let $U$ be the affine cubic surface over $\Qbb$ defined by
$$ x^{2} + y^{2} + z^{2} + xyz = ax + by + cz + d, $$
where $a,b,c,d$ are defined by $(2)$. Assume that its natural compactification $X \subset \Pbb^{3}$ is smooth over $\Qbb$. Set $X_{E} := X \times_{\Qbb} E$ and $U_{E} := U \times_{\Qbb} E$. Then the natural map $\textup{Br}\,X_{E} \rightarrow \textup{Br}\,U_{E}$ is an isomorphism. Moreover, $U$ has trivial transcendental Brauer group over $\Qbb$.
\end{proposition}

\begin{proof}
The proof is inspired by that of \cite[Proposition 4.1]{LM20}. Let $\mathcal{B} \in \textup{Br}\,U_{E}$ be a non-constant Brauer element. Multiplying $\mathcal{B}$ by a constant algebra if necessary, by Proposition 3.3 and Theorem 3.8, we may assume that $\mathcal{B}$ has order dividing $4$ (note that under our assumption of $k$, the field extension $E$ is totally real and thus contains no nontrivial roots of unity). In order to show that $\mathcal{B} \in \textup{Br}\,X_{E}$, we only need to show that $\mathcal{B}$ is unramified along the three lines $L_{i}$ on $X$ by Grothendieck's purity theorem (see \cite[Theorem 3.7.2]{CTS21}).

Let $L = L_{1}$ and $C = L_{2} \cup L_{3}$. Let $L' = L\,\backslash\,C$. Note that $L$ meets $C$ at two rational points, so $L'$ is non-canonically isomorphic to $\Gmbb$. Let the point $(x:y:z:t) = (0:1:1:0) \in L'$ be the identity element of the group law. Then an isomorphism with $\Gmbb$ is realized via the following homomorphism: 
\begin{equation}
\Gmbb \rightarrow X, \hspace{1cm} u \mapsto (0:u:1:0).
\end{equation}
The residue of $\mathcal{B}$ along $L$ lies inside $\textup{H}^{1}(L',\Zbb/2\Zbb)$. Assume by contradiction that the residue is nontrivial. Since the order of $\mathcal{B}$ is a power of $2$ dividing $4$, then we can assume that the residue has order $2$ (up to replacing $\mathcal{B}$ by $2\mathcal{B}$). This means that the residue corresponds to some irreducible degree $2$ finite étale cover $f : L'' \rightarrow L'$. 

Over the field $E$, the conic fiber $C$ over the coordinate $\left( x = \ds\frac{k_{1}k_{2} + \sqrt{(k_{1}^{2}-4)(k_{2}^{2}-4)}}{2} t \right)$ is split, i.e. a union of two lines over $E$. These lines meet $L$ at the points 
$$ Q_{\pm} = (0 : \ds\frac{(k_{1} \pm \sqrt{k_{1}^{2}-4})(k_{2} \pm \sqrt{k_{2}^{2}-4})}{4} : 1 : 0). $$
Let $C_{+}$ be the irreducible component of $C$ containing $Q_{+}$, i.e., $C_{+} = \ell_{1}(1,1)$. Consider the restriction of $\mathcal{B}$ to $C_{+}$. This is well-defined outside of $Q_{+}$, and since $C_{+}\,\backslash\,Q_{+} \simeq \Abb^{1}$ has constant Brauer group, $\mathcal{B}$ actually extends to all of $C_{+}$. As $C_{+}$ meets $L$ transversely, by the functoriality of residues (\cite[Section 3.7]{CTS21}) we deduce that the residue of $\mathcal{B}$ at $Q_{+}$ is also trivial, so the fiber $f^{-1}(Q_{+})$ consists of exactly two rational points. This implies that $L''$ is geometrically irreducible, hence $L'' \cong \Gmbb$ non-canonically.

Now by choosing a rational point over $Q_{+}$ and using the above group homomorphism, we may therefore identify the degree $2$ cover $L'' \rightarrow L'$ with the map
\begin{equation}
\Gmbb \rightarrow X, \hspace{1cm} u \mapsto (0:u^{2}:1:0).
\end{equation}
However, our assumptions on $k$ imply that $\ds\frac{(k_{1} + \sqrt{k_{1}^{2}-4})(k_{2} + \sqrt{k_{2}^{2}-4})}{4}$ is not a square in $E^{\times}$, which gives a contradiction. Thus the residue of $\mathcal{B}$ along $L$ is trivial, and the same holds for the other lines. We conclude that $\mathcal{B}$ is everywhere unramified, hence $\mathcal{B} \in \textup{Br}\,X_{E}$. 
\\~\\
\indent Now let $B \in \textup{Br}\,U$ be a non-constant element. Then over the field extension $E$, the corresponding image of $B$ comes from $\textup{Br}\,X_{E}$ by the above argument. As $\textup{Br}\,\overline{X} = 0$, it is clear that $B$ is algebraic. The result follows.
\end{proof}

\begin{remark}
Note that $(\textup{Br}\,U_{\overline{\Qbb}})^{\textup{Gal}(\overline{\Qbb}/\Qbb)} = \Zbb/2\Zbb$ by Proposition 3.3. However, in the above proposition, the Galois invariant element of order $2$ does not descend to a Brauer group element over $\Qbb$, which is also the case in \cite[Proposition 4.1]{LM20}.
\end{remark}

\section{The Brauer--Manin obstruction}
\subsection{Review of the Brauer--Manin obstruction}
Here we briefly recall how the Brauer--Manin obstruction works in our setting, following \cite[Section 8.2]{Poo17} and \cite[Section 1]{CTX09}. For each place $v$ of $\Qbb$ there is a pairing
$$ U(\Qbb_{v}) \times \textup{Br}\,U \rightarrow \Qbb/\Zbb $$
coming from the local invariant map 
$$ \textup{inv}_{v} : \textup{Br}\,\Qbb_{v} \rightarrow \Qbb/\Zbb $$ from local class field theory (this is an isomorphism if $v$ is a prime number). This pairing is locally constant on the left by \cite[Proposition 8.2.9]{Poo17}.
For \textbf{integral points}, any element $\alpha \in \textup{Br}\,U$ pairs trivially on $\mathcal{U}(\Zpbb)$ for almost all primes $p$, so we obtain a pairing $U(\textbf{\textup{A}}_{\Qbb}) \times \textup{Br}\,U \rightarrow \Qbb/\Zbb$. As the local pairings are locally constant, we obtain a well-defined pairing 
$$ \mathcal{U}(\textbf{\textup{A}}_{\Zbb})_{\bullet} \times \textup{Br}\,U \rightarrow \Qbb/\Zbb. $$
For $B \subseteq \textup{Br}\,U$, let $\mathcal{U}(\textbf{\textup{A}}_{\Zbb})_{\bullet}^{B}$ be the left kernel with respect to $B$, and let $\mathcal{U}(\textbf{\textup{A}}_{\Zbb})_{\bullet}^{\textup{Br}} = \mathcal{U}(\textbf{\textup{A}}_{\Zbb})_{\bullet}^{\textup{Br}\,U}$. By abuse of notation, from now on we write the \emph{reduced} Brauer--Manin set $\mathcal{U}(\textbf{\textup{A}}_{\Zbb})_{\bullet}^{B}$ in the standard way as $\mathcal{U}(\textbf{\textup{A}}_{\Zbb})^{B}$. Note that the set $\mathcal{U}(\textbf{\textup{A}}_{\Zbb})^{B}$ depends only on the image of $B$ in the quotient $\textup{Br}\,U/\textup{Br}_{0}\,U$. By Theorem 3.8, the map $\langle \mathcal{A} \rangle \rightarrow \textup{Br}_{1}\,U/\textup{Br}\,\Qbb$ is an isomorphism, hence $\mathcal{U}(\textbf{\textup{A}}_{\Zbb})^{\textup{Br}_{1}} = \mathcal{U}(\textbf{\textup{A}}_{\Zbb})^{\textup{Br}_{1}\,U} = \mathcal{U}(\textbf{\textup{A}}_{\Zbb})^{\mathcal{A}}$. We have the inclusions $\mathcal{U}(\Zbb) \subseteq \mathcal{U}(\textbf{\textup{A}}_{\Zbb})^{\mathcal{A}} \subseteq \mathcal{U}(\textbf{\textup{A}}_{\Zbb})$, so that $\mathcal{A}$ can obstruct the \textit{integral Hasse principle} or \textit{strong approximation} on $\mathcal{U}$.

Let $V$ be dense Zariski open in $U$. As $U$ is smooth, the set $V(\Qbb_{v})$ is dense in $U(\Qbb_{v})$ for all places $v$. Moreover, $\mathcal{U}(\Zpbb)$ is open in $U(\Qpbb)$, hence $V(\Qpbb) \cap \mathcal{U}(\Zpbb)$ is dense in $\mathcal{U}(\Zpbb)$. As the local pairings are locally constant, we may restrict our attention to $V$ to calculate the local invariants of a given element in $\textup{Br}\,U$. In particular, here we take the open subset $V$ given by $((k_{1}^{2}-4)(k_{2}^{2}-4) - (2x-k_{1}k_{2})^{2})((k_{1}^{2}-4)(k_{4}^{2}-4) - (2y-k_{1}k_{4})^{2})((k_{1}^{2}-4)(k_{3}^{2}-4) - (2z-k_{1}k_{3})^{2}) \not= 0$.

\subsection{Brauer--Manin obstruction from a quaternion algebra}
Now we consider the Markoff-type cubic surfaces $U$ over $\Qbb$ and their integral models $\mathcal{U}$ over $\Zbb$ defined by the equation $(1)$:
$$ x^{2} + y^{2} + z^{2} + xyz = ax + by + cz + d, $$
where $a,b,c,d \in \Zbb$ are defined as in $(2)$ with $k \in \Zbb^{4}$. Set $f := x^{2} + y^{2} + z^{2} + xyz - ax - by - cz - d \in \Zbb[x,y,z]$. First of all, we study the existence of local integral points on those affine cubic surfaces given by $f=0$. 

\begin{proposition}[Assumption A]
If $k = (k_{1},k_{2},k_{3},k_{4}) \in (\Zbb\,\backslash\,[-2,2])^{4}$ satisfies $k_{1} \equiv -1$ \textup{mod} $16$, $k_{i} \equiv 5$ \textup{mod} $16$ and $k_{1} \equiv 1$ \textup{mod} $9$, $k_{i} \equiv 5$ \textup{mod} $9$ for $2 \leq i \leq 4$, such that $(k_{i},k_{j}) = 1$, $(k_{i}^{2}-4,k_{j}^{2}-4) = 3$ for $1 \leq i \not= j \leq 4$, and $(k_{1}^{2}-2,k_{2}^{2}-2,k_{3}^2-2,k_{4}^{2}-2) = 1$, then $\mathcal{U}(\textbf{\textup{A}}_{\Zbb}) \not= \emptyset$.
\end{proposition}

\begin{proof}
With our specific choice of $k$ in the assumption, we obtain:
\begin{enumerate}
	\item[(1)] Prime powers of $p = 2$: The only solution modulo $2$ is the singular $(1,0,0)$ (up to permutation). However, we find a solution $(1,0,0)$ modulo $8$ with $2x + yz \equiv 2$ mod $8$, so twice the order valuation at $2$ of the partial derivative at $x$ is less than the order valuation at $2$ of $f(x,y,z)$. This solution then lifts to a $2$-adic integer solution by Hensel's lemma (fixing the variables $y,z$).
	
	\item[(2)] Prime powers of $p = 3$: We find the non-singular solution $(1,0,0)$, which lifts to a $3$-adic integer solution by Hensel's lemma (fixing the variables $y,z$).
	
	\item[(3)] Prime powers of $p \geq 5$: We would like to find a non-singular solution modulo $p$ of the equations $f = 0$ which does not satisfy simultaneously 
	$$ df=0 : 2x+yz=a, 2y+xz=b, 2z+xy=c. $$ 
	For simplicity, we will find a sufficient condition for the existence of a non-singular solution whose at least one coordinate is zero. First, it is clear that the equation $f=0$ always has a solution whose one coordinate is $0$: indeed, take $z=0$, then $f=0$ is equivalent to $(2x-a)^{2} + (2y-b)^{2} = (k_{1}^{2}+k_{3}^{2}-4)(k_{2}^{2}+k_{4}^{2}-4)$, and every element in $\Fpbb$ can be expressed as a sum of two squares. Now assume that all such points in $U(\Fpbb)$ with one coordinate equal to $0$ are singular. Then modulo $p$, if $z=0$ the required equations become $f=0, (2x=a, 2y=b, xy = c)$, plus all the permutations for $x=0$ and $y=0$. For convenience, we drop the phrase ``modulo $p$''. From these equations, we get
\begin{equation}	
	(k_{1}^{2}+k_{3}^{2}-4)(k_{2}^{2}+k_{4}^{2}-4) = 0, ab=4c\,(=4xy)
\end{equation}	
plus all the permutations, respectively 
\begin{equation}	
	(k_{1}^{2}+k_{2}^{2}-4)(k_{3}^{2}+k_{4}^{2}-4) = 0, bc=4a
\end{equation}
and 
\begin{equation}	
	(k_{1}^{2}+k_{4}^{2}-4)(k_{2}^{2}+k_{3}^{2}-4) = 0, ac=4b. 
\end{equation}
We will choose $k$ which does not satisfy all of these equations simultaneously to get a contradiction to our assumption at the beginning.
	\\~\\
	Now if $k$ satisfies all the above equations, we first require that $(k_{i}^{2}-4,k_{j}^{2}-4) = 3$ for any $1 \leq i \not= j \leq 4$. Next, without loss of generality (WLOG), assume in $(22)$ that $k_{1}^{2}+k_{4}^{2}-4=0$ with $ac = 4b = (k_{1}^{2}+k_{4}^{2})b$, then putting the formulas of $a,b,c$ in $ac - (k_{1}^{2}+k_{4}^{2})b = 0$ gives us $k_{1}k_{4}(k_{2}^{2}+k_{3}^{2}-k_{1}^{2}-k_{4}^{2}) = 0$, hence either $k_{2}^{2}+k_{3}^{2}-4=0$ or $k_{1}k_{4}=0$. (We also have similar equations for all the other cases.)
	
	If $k_{2}^{2}+k_{3}^{2}-4=0$, then $\sum_{i=1}^{4} k_{i}^{2} - 8 = 0$, and so from $(20),(21)$ we obtain $k_{i}^{2}+k_{j}^{2}-4=0$ for any $i \not= j$, hence $k_{i}^{2} - 2 = 0$ for all $1 \leq i \leq 4$. Otherwise, if $k_{2}^{2}+k_{3}^{2}-4 \not= 0$ but $k_{1}k_{4}=0$, WLOG assume in $(22)$ that $k_{1} = 0$, then $k_{4}^{2}-4 = 0$ and so from $(20),(21)$ we have the following possibilities: $k_{2}^{2} - 4 = k_{3}^{2} - 4 = 0$, or $k_{2} = k_{3} = 0$. Therefore, we immediately deduce a sufficient condition for the nonexistence of singular solutions modulo $p \geq 5$: $(k_{i},k_{j}) = 1$, $(k_{i}^{2}-4,k_{j}^{2}-4) = 3$ for any $1 \leq i \not= j \leq 4$, and $(k_{1}^{2}-2,k_{2}^{2}-2,k_{3}^{2}-2,k_{4}^{2}-2) = 1$.
	
	As a result, assuming the hypothesis of the proposition, it is clear that $U(\Fpbb)$ has a smooth point, which then lifts to a $p$-adic integral point by Hensel's lemma (with respect to the variable at which the partial derivative is nonzero modulo $p$, fixing other variables).
\end{enumerate}
\end{proof} 

We keep Assumption A to ensure that $\mathcal{U}(\textbf{A}_{\Zbb}) \not= \emptyset$ and study the Brauer--Manin obstruction to the existence of integral points. Note that our specific choice of $k$ implies that $[E:\Qbb]=16$. Indeed, it is clear that with our assumption, $k_{i}^{2} - 4 > 0$ is not a square in $\Qbb$ for $1 \leq i \leq 4$. Then $[E:\Qbb] < 16$ if and only if there are some $i \not= j$ such that $(k_{i}^{2}-4)(k_{j}^{2}-4)$ is a square in $\Qbb$; however, since $(k_{i}^{2}-4,k_{j}^{2}-4) = 3$ and $k_{i}^{2} - 4 = 3l^{2}$ does not have any solution modulo $8$ with $k_{i}$ odd, that cannot happen. 

Now let us calculate the local invariants of the following quaternion algebras as elements of the algebraic Brauer group $\textup{Br}_{1}\,U$:

\begin{equation*}
	\begin{cases*}
	\mathcal{A}_{1} = \textup{Cor}^{F_{1}}_{\Qbb} \left( x - \frac{k_{1}k_{2} + \sqrt{(k_{1}^{2}-4)(k_{2}^{2}-4)}}{2}, (k_{1}\sqrt{k_{2}^{2}-4} + k_{2}\sqrt{k_{1}^{2}-4})^{2} \right), \\
	\mathcal{A}_{2} = \textup{Cor}^{F_{3}}_{\Qbb} \left( y - \frac{k_{1}k_{4} + \sqrt{(k_{1}^{2}-4)(k_{4}^{4}-4)}}{2}, (k_{1}\sqrt{k_{4}^{2}-4} + k_{4}\sqrt{k_{1}^{2}-4})^{2} \right), \\
	\mathcal{A}_{3} = \textup{Cor}^{F_{2}}_{\Qbb} \left( z - \frac{k_{1}k_{3} + \sqrt{(k_{1}^{2}-4)(k_{3}^{2}-4)}}{2}, (k_{1}\sqrt{k_{3}^{2}-4} + k_{3}\sqrt{k_{1}^{2}-4})^{2} \right),
	\end{cases*} 
\end{equation*} 
where $F_{i} = \Qbb(\sqrt{(k_{1}^{2}-4)(k_{i+1}^{2}-4)})$ for $1 \leq i \leq 3$. Now for each $i$, we have the local invariant map
$$ \textup{inv}_{p}\,\mathcal{A}_{i} : U(\Qpbb) \rightarrow \Zbb/2\Zbb, \hspace{1cm} x \mapsto \textup{inv}_{p}\,\mathcal{A}_{i}(x). $$

\begin{lemma}[Assumption B]
Let $k = (k_{1},k_{2},k_{3},k_{4}) \in \Zbb^{4}$ and $p \geq 11$ be a prime such that $(2,p)_{p} = 1/2$ and $k_{1} - 2 \equiv p$ \textup{mod} $p^{3}$. Then there exist $k_{2},k_{3},k_{4}$ such that $k_{2}^{2} - 4$ is a quadratic nonresidue \textup{mod} $p$, $k_{3} \equiv k_{4} \not\equiv 0, 2$ \textup{mod} $p$, $k_{2} - 2 \not\equiv 2k_{3}$ \textup{mod} $p$, and $k_{2},k_{3},k_{4}$ satisfy the cubic equation defining an affine Markoff surface $S$:
$$ k_{2}^{2} + k_{3}^{2} + k_{4}^{2} - k_{2}k_{3}k_{4} \equiv (p + 2)^{2} \hspace{0.25cm} \textup{mod}\, p^{2}. $$ 
Furthermore, for any $k$ satisfying all the above conditions (Assumption B) and Assumption A, the local invariant map of the quaternion algebra $\mathcal{A}_{1}$ at $p$ is surjective. 
\end{lemma}

\begin{proof}
For convenience (and by abuse of notation), we will also write $\mathcal{A}_{1}$ (resp. $F_{1}$) as $\mathcal{A}$ (resp. $F$). We denote $(k_{1}\sqrt{k_{2}^{2}-4} + k_{2}\sqrt{k_{1}^{2}-4})$ by $\alpha_{1}$ and write it simply as $\alpha$. Set $D := (k_{1}^{2}-4)(k_{2}^{2}-4)$.

First of all, over $F = \Qbb(\sqrt{D})$ the affine cubic equation of $\mathcal{U}$ can be rewritten equivalently as 
\begin{equation}
\begin{aligned}
f(x,y,z,k_{1},k_{2},k_{3},k_{4}) = \\
\left( x - \frac{k_{1}k_{2} + \sqrt{D}}{2} \right) \left( x - \frac{k_{1}k_{2} - \sqrt{D}}{2} - k_{3}k_{4} + yz \right) &+ \left( y + \frac{k_{1}k_{2} + \sqrt{D}}{4}z - \frac{b}{2} \right)^{2} \\ &- \frac{\alpha^{2}}{16}\left( z - \frac{2b(k_{1}k_{2}+\sqrt{D}) - 8c}{\alpha^{2}} \right)^{2} = 0,
\end{aligned}
\end{equation}
for all $(k_{1}, k_{2}, k_{3}, k_{4}) \in \Zbb^{4}$ satisfying our hypothesis. Since $\textup{v}_{p}(D) = 1$, $p$ ramifies over $F = \Qbb(\sqrt{D})$, i.e., $p\mathcal{O}_{F} = \mathfrak{p}^{2}$ where $\mathfrak{p}$ is a nonzero prime ideal of $\mathcal{O}_{F}$. Therefore, we have $\textup{N}\mathfrak{p} = p$ and 
$$ \mathcal{O}_{F} / \mathfrak{p} \cong \Fpbb = \Zbb/p\Zbb. $$
Following the proof of \cite[Proposition 5.5]{LM20}, for all $u_{2},u_{3},u_{4} \in \Fpbb^{*}$, there exists an $\Fpbb$-point on the variety
$$ (2k_{2}-k_{3}k_{4})^{2} = (k_{3}^{2}-4)(k_{4}^{2}-4), k_{2}-2=u_{2}v_{2}^{2}, k_{3}-2=u_{3}v_{3}^{2}, k_{4}-2 = u_{4}v_{4}^{2} $$
which satisfy $v_{2}v_{3}v_{4} \not= 0$. As $(k_{2}-2)(k_{3}-2)(k_{4}-2) \not= 0$, we know from \cite[Lemma 5.3]{LM20} that this gives rise to a smooth $\Fpbb$-point of $S$, hence a $\Zbb/p^{2}\Zbb$-point with the same residue modulo $p$ by Hensel's lemma. Now to construct the given $\Fpbb$-point, we restrict our attention to the subvariety given by $k_{3} = k_{4}$ mod $p$ then assume that $u_{3} = u_{4}$ mod $p$ and $u_{2}$ is a quadratic nonresidue mod $p$. The above equations then become
$$ (2k_{2}-k_{3}k_{4})^{2} = (k_{3}^{2}-4)^{2}, k_{2}-2 = u_{2}v_{2}^{2}, k_{3}-2 = u_{3}v_{3}^{2}. $$
Factoring the left hand side, it suffices to solve the equations
$$ k_{2} = k_{3}^{2}-2, k_{2}-2 = u_{2}v_{2}^{2}, k_{3}-2 = u_{3}v_{3}^{2}. $$
This then gives the equation of an affine curve
$$ u_{2}v_{2}^{2} = u_{3}^{2}v_{3}^{4} + 4u_{3}v_{3}^{2}. $$
By the argument in the proof of \cite[Proposition 5.5]{LM20}, this affine curve has a unique singular point $(v_{2},v_{3}) = (0,0)$ and has $p-2$ many $\Fpbb$-points. Of these points at most $3$ satisfy $v_{2}v_{3} = 0$, and $k_{3} = 0$ gives only at most $4$ points, hence providing $p - 2 - 3 - 4 = p - 9 > 0$, there exists an $\Fpbb$-point $(k_{2},k_{3},k_{4})$ with the properties:
$$ (k_{2}-2)(k_{3}-2)(k_{4}-2) \not= 0, k_{3} = k_{4} \not= 0. $$
Since $k_{3} \not= 0$, we have $k_{2} = k_{3}^{2}-2 \not= -2$, so $k_{2}^{2}-4 \not= 0$ in $\Fpbb$. Moreover, as $u_{2}$ is a quadratic nonresidue modulo $p$, so is $k_{2}-2 = k_{3}^{2}-4$. Since $k_{2}+2 = k_{3}^{2}$ is a nonzero square in $\Fpbb$, we deduce that $k_{2}^{2}-4$ is a quadratic nonresidue modulo $p$, as required.

Now after lifting from the smooth $\Fpbb$-point to a $\Zbb/p^{2}\Zbb$-point (with the same residue modulo $p$), if $k_{2}-2 \equiv 2k_{3}$ mod $p$ then we can take $k_{2}' = k_{2}$, $k_{3}' = -k_{3}$ and $k_{4}' = -k_{4}$ (which still satisfy the cubic equation mod $p^{2}$) to have $k'_{2} - 2 \not\equiv 2k_{3}'$ mod $p$ since $k_{3} \not\equiv 0$ mod $p$ and $p$ is odd. Therefore, there exist integers $k_{2}, k_{3}, k_{4}$ satisfying the required properties by the Chinese Remainder Theorem.
\\~\\
\indent Now assume that $k = (k_{1},k_{2},k_{3},k_{4})$ satisfies all the above congruence conditions and Assumption A. Since $p\mathcal{O}_{F} = \mathfrak{p}^{2}$ and $\mathcal{O}_{F}/\mathfrak{p} \cong \Zbb/p\Zbb$, by Hensel's lemma we deduce that when $k_{2}^{2}-4$ is a quadratic nonresidue modulo $p$ over $\Zbb$, then $k_{2}^{2}-4$ is a quadratic nonresidue modulo $\mathfrak{p}$ over $\mathcal{O}_{F}$. Recalling that $k_{1} - 2 \equiv p$ mod $p^{3}$, we can choose a $\Zbb/p^{3}\Zbb$-point $(x,y,z)$ of $\mathcal{U}$ such that 
$$ x - k_{2} = p, y - k_{4} = p, z - k_{3} = p, $$
from which we have (in $\Zbb/p^{3}\Zbb$):
$$
\begin{aligned} 
&f(x,y,z,k_{1},k_{2},k_{3},k_{4}) = f(x,y,z,2,k_{2},k_{3},k_{4}) - \sum_{i \in \{2,3,4\}} pk_{i}(k_{i}+p) + 4p+p^{2} + pk_{2}k_{3}k_{4}\\
&= (x-k_{2})(x-k_{2}+yz-k_{3}k_{4}) + \left(y+\frac{k_{2}z}{2}-\frac{b}{2} \right)^{2} - \frac{k_{2}^{2}-4}{4}(z-k_{3})^{2} \\& - \sum_{i \in \{2,3,4\}} pk_{i}(k_{i}+p) + 4p+p^{2} + pk_{2}k_{3}k_{4}\\
&= p(p+p^{2}+p(k_{3}+k_{4})) + p^{2}\left(1+\frac{k_{2}}{2}\right)^{2} - p^{2}\left(\frac{k_{2}^{2}-4}{4}\right) - \sum_{i \in \{2,3,4\}} pk_{i}(k_{i}+p) + 4p+p^{2} + pk_{2}k_{3}k_{4}\\
&= p(p+p^{2}+p(k_{3}+k_{4})) + p(2p+pk_{2}) - \sum_{i \in \{2,3,4\}} p(k_{i}^{2}+pk_{i}) + p(4+p+k_{2}k_{3}k_{4})\\
&= p(4+4p+p^{2} - (k_{2}^{2}+k_{3}^{2}+k_{4}^{2}) + k_{2}k_{3}k_{4})\\
&= p((p+2)^{2} - (k_{2}^{2}+k_{3}^{2}+k_{4}^{2} - k_{2}k_{3}k_{4})) = 0,
\end{aligned} 
$$ 
and
$$ f'_{x}(x,y,z) = 2x + yz - k_{1}k_{2} - k_{3}k_{4} = 2(x-k_{2}) - pk_{2} + p^{2} + p(k_{3}+k_{4}) = p(p+2-k_{2}+k_{3}+k_{4}). $$ 
Therefore $\textup{v}_{p}(f) \geq 3 > 2 = 2\textup{v}_{p}(f'_{x})$, by Hensel's lemma (fixing the variables $y,z$), this point lifts to a $\Zpbb$-point $(x,y,z)$ (abuse of notation) with the same residue modulo $p^{2}$. Then since $\textup{v}_{\mathfrak{p}}(\sqrt{D}) = \textup{v}_{p}(D) = 1$, the local invariant at $p$ of $\mathcal{A}$ at this point is equal to $$ \left( x - \frac{k_{1}k_{2} + \sqrt{D}}{2}, \alpha^{2} \right)_{\mathfrak{p}} = \left( x - \frac{k_{1}k_{2} + \sqrt{D}}{2}, k_{2}^{2}-4 \right)_{\mathfrak{p}} = 1/2 $$ by the formulae in \cite[Proposition II.1.4 and Proposition III.3.3]{Neu13}.

It is also clear that there exists a $\Zpbb$-point such that $x-k_{2}$ is not divisible by $p$, which gives the local invariant of $\mathcal{A}$ at $p$ the value $0$. Indeed, if every $\Zpbb$-point satisfies that $x-k_{2}$ is divisible by $p$, then as $k_{1}-2 \equiv 0$ mod $p$, the affine equation of $\mathcal{U}$ over $\Fpbb$ becomes
$$
\begin{aligned} 
(x-k_{2})(x-k_{2}+yz-k_{3}k_{4}) + \left(y+\frac{k_{2}z}{2}-\frac{b}{2} \right)^{2} - \frac{k_{2}^{2}-4}{4}(z-k_{3})^{2} = 0.
\end{aligned} 
$$ 
From the fact that $k_{2}^{2}-4$ is a quadratic nonresidue mod $p$, one must obtain that (in $\Fpbb$):
$$ y + \frac{k_{2}z}{2} - \frac{b}{2} = 0, z - k_{3} = 0, $$
which gives $y = k_{4}, z = k_{3}$. Therefore, in $\Fpbb$ we have
$$ 2x+yz = k_{1}k_{2}+k_{3}k_{4}, 2y+xz = k_{1}k_{4}+k_{2}k_{3}, 2z+xy = k_{1}k_{3}+k_{2}k_{4}, $$ which implies that $f'_{x} = f'_{y} = f'_{z} = 0$. This cannot be true for every $\Zpbb$-point $(x,y,z)$ of $\mathcal{U}$, since by Assumption A we have proved in Proposition 4.1 that there always exists a non-singular solution modulo $p$ of the affine Markoff-type cubic equation (with at least one coordinate zero) which lifts to a $\Zpbb$-point by Hensel's lemma.

In conclusion, the local invariant map of $\mathcal{A}$ at the prime $p$ satisfying our hypothesis is indeed surjective.
\end{proof}

\begin{theorem}
Let $k \in (\Zbb\,\backslash\,[-2,2])^{4}$ satisfy the hypotheses of Assumptions A and B. Then we have a Brauer--Manin obstruction to strong approximation for $\mathcal{U}$ given by the class of $\mathcal{A} = \mathcal{A}_{1}$ in $\textup{Br}_{1}\,U/\textup{Br}\,\Qbb$ (i.e. $\mathcal{U}(\textbf{\textup{A}}_{\Zbb}) \not= \mathcal{U}(\textbf{\textup{A}}_{\Zbb})^{\mathcal{A}} = \mathcal{U}(\textbf{\textup{A}}_{\Zbb})^{\textup{Br}_{1}}$) and no \emph{algebraic} Brauer--Manin obstruction to the integral Hasse principle for $\mathcal{U}$ (i.e. $\mathcal{U}(\textbf{\textup{A}}_{\Zbb})^{\textup{Br}_{1}} \not= \emptyset$).
\end{theorem}

\begin{proof}
By abuse of notation, we also denote the class of $\mathcal{A}$ by the Brauer group element itself. For any point $\textbf{u} = (x,y,z) \in \mathcal{U}(\textbf{A}_{\Zbb})$, from the above Lemma we can find a local point $\textbf{u}_{p} = (x_{p},y_{p},z_{p}) \in \mathcal{U}(\Zbb_{p})$ such that $\textup{inv}_{p}\,\mathcal{A}(\textbf{u}_{p}) = - \sum_{q \not= p} \textup{inv}_{q} \mathcal{A}(\textbf{u})$ and $\textbf{u}_{p}' = (x_{p}',y_{p}',z_{p}') \in \mathcal{U}(\Zbb_{p})$ such that $\textup{inv}_{p}\,\mathcal{A}(\textbf{u}_{p}') = 1/2 - \sum_{q \not= p} \textup{inv}_{q}\,\mathcal{A}(\textbf{u})$. Then we obtain a point in $\mathcal{U}(\textbf{A}_{\Zbb})^\mathcal{A}$ by replacing the $p$-part of $\textbf{u}$ by $\textbf{u}_{p}$, and a point in $\mathcal{U}(\textbf{A}_{\Zbb})$ but not in $\mathcal{U}(\textbf{A}_{\Zbb})^\mathcal{A}$ by replacing the $p$-part of $\textbf{u}$ by $\textbf{u}_{p}'$.

Therefore, we have a Brauer--Manin obstruction to strong approximation and no algebraic Brauer--Manin obstruction to the integral Hasse principle for $\mathcal{U}$.
\end{proof}

\begin{example}
For $p = 11$: Take $k_{1} = 11^3.102+13 = 135775, k_{2} = 11^2.112+111 = 13663, k_{3} = 11^2.119+6 = 14405, k_{4} = 11^2.(119+144)+6 = 31829$.
\end{example}

\begin{remark}
In addition, if $k$ satisfies Assumption 3.3 from the previous section, then we can drop the term ``algebraic'' from the statement of Theorem 4.3, and we will have no Brauer--Manin obstruction to the integral Hasse principle.
\end{remark}

\subsection{Some counting results}
In this part, we compute the number of examples of existence for local integral points as well as the number of counterexamples to strong approximation for the Markoff-type cubic surfaces in question which can be explained by the Brauer--Manin obstruction. More precisely, we consider the natural density of $k \in (\Zbb\,\backslash\,[-2,2])^{4}$ satisfying $k_{1} \equiv -1$ mod $16$, $k_{i} \equiv 5$ mod $16$ and $k_{1} \equiv 1$ mod $9$, $k_{i} \equiv 5$ mod $9$ for $2 \leq i \leq 4$, such that $(k_{i},k_{j}) = 1$, $(k_{i}^{2}-4,k_{j}^{2}-4) = 3$ for $1 \leq i \not= j \leq 4$, and $(k_{1}^{2}-2,k_{2}^{2}-2,k_{3}^{2}-2,k_{4}^{2}-2) = 1$ (then satisfying the additional hypothesis in Lemma 4.2). Note that the finite number of $k \in [-2,2]^{4}$ is negligible here.

In fact, for now we can only give an asymptotic lower bound. To get a lower bound, it is enough to count the number of examples satisfying stronger conditions, namely the congruence conditions for the $k_{i}$ and the common divisor condition 
$$ \textup{gcd}(k_{i}(k_{i}^{2}-2)(k_{i}^{2}-4),k_{j}(k_{j}^{2}-2)(k_{j}^{2}-4)) = 3 $$ for $1 \leq i \not= j \leq 4$. To do this, we make use of a natural generalization of Ekedahl--Poonen's formula in \cite[Theorem 3.8]{Poo03} as follows.

\begin{proposition}
Let $f_{i} \in \Zbb[x_{1},\dots,x_{n}]$, $1 \leq i \leq s$ for some $s \in \Zbb_{>1}$, be $s$ polynomials that are mutually relatively prime as elements of $\Qbb[x_{1},\dots,x_{n}]$. For each $1 \leq i \not= j \leq s$, let
$$ \mathcal{R}_{i,j} := \{ a \in \Zbb^{n} : \textup{gcd}(f_{i}(a),f_{j}(a))=1 \}. $$
Then $\mu(\bigcap_{i \not= j} \mathcal{R}_{i,j}) = \prod_{p} (1 - c_{p}/p^{n})$, where $p$ ranges over all primes of $\Zbb$, and $c_{p}$ is the number of $x \in (\Zbb/p\Zbb)^{n}$ satisfying at least one of $f_{i}(x) = f_{j}(x) = 0$ in $\Zbb/p\Zbb$ for $1 \leq i\not=j \leq s$. 
\end{proposition}

\begin{proof}
We also have a generalization of \cite[Lemma 5.1]{Poo03} as follows.

\begin{lemma}
Let $f_{i} \in \Zbb[x_{1},\dots,x_{n}]$, $1 \leq i \leq s$ for some $s \in \Zbb_{>1}$, be $s$ polynomials that are mutually relatively prime as elements of $\Qbb[x_{1},\dots,x_{n}]$. For each $1 \leq i \not= j \leq s$, let
$$ \mathcal{Q}_{i,j,M} := \{ a \in \Zbb^{n} : \exists\,p\;\textup{such that}\; p \geq M\;\textup{and}\; p\,|\,f_{i}(a),f_{j}(a) \}. $$
Then $\lim_{M \rightarrow \infty} \overline{\mu}(\bigcup_{i \not= j} \mathcal{Q}_{i,j,M}) = 0$.
\end{lemma}

\begin{proof}
The result immediately follows from the inequality $\overline{\mu}(\bigcup_{i \not= j} \mathcal{Q}_{i,j,M}) \leq \sum_{i \not= j} \overline{\mu}(\mathcal{Q}_{i,j,M})$ and the original result of \cite[Lemma 5.1]{Poo03}.
\end{proof}

Next, we proceed similarly as in the proof of \cite[Theorem 3.1]{Poo03}. Let $P_{M}$ denote the set of prime numbers of $\Zbb$ such that $p < M$. Approximate $\mathcal{R}_{i,j}$ by
$$ \mathcal{R}_{i,j,M} := \{ a \in \Zbb^{n} : f_{i}(a)\;\textup{and}\;f_{j}(a)\;\textup{are not both divisible by any prime}\; p \in P_{M} \}. $$
Define the ideal $I$ as the product of all $(p)$ for $p \in P_{M}$. Then $\mathcal{R}_{i,j,M}$ is a union of cosets of the subgroup $I^{n} \subset \Zbb^{n}$. Hence $\mu(\bigcap_{i \not= j} \mathcal{R}_{i,j,M})$ is the fraction of residue classes in $(\Zbb/I)^{n}$ in which for all $p \in P_{M}$, for all $1 \leq i\not=j \leq s$, at least one of $f_{i}(a)$ and $f_{j}(a)$ is nonzero modulo $p$. Applying the Chinese Remainder Theorem, we obtain that $\mu(\bigcap_{i \not= j} \mathcal{R}_{i,j,M}) = \prod_{p \in P_{M}} (1 - c_{p}/p^{n})$. By the above lemma,
\begin{equation}
\mu(\bigcap_{i \not= j} \mathcal{R}_{i,j}) = \lim_{M \rightarrow \infty} \mu(\bigcap_{i \not= j} \mathcal{R}_{i,j,M}) = \prod_{p} (1 - c_{p}/p^{n}).
\end{equation}
Since $f_{i}$ and $f_{j}$ are relatively prime as elements of $\Qbb[x_{1},\dots,x_{n}]$ for any $i \not= j$, there exists a nonzero $u \in \Zbb$ such that every $f_{i} = f_{j} = 0$ defines a subscheme of $\Abb^{n}_{\Zbb[1/u]}$ of codimension at least $2$. Thus $c_{p} = O(p^{n-2})$ as $p \rightarrow \infty$, and the product converges.
\end{proof}

\begin{theorem}
Let $\mathcal{U}$ be the affine scheme over $\Zbb$ defined by 
\begin{equation*}
x^{2} + y^{2} + z^{2} + xyz = ax + by + cz + d,
\end{equation*}
where 
\begin{equation*}
	\begin{cases*}
	a = k_{1}k_{2} + k_{3}k_{4} \\
	b = k_{1}k_{4} + k_{2}k_{3} \\
	c = k_{1}k_{3} + k_{2}k_{4}
	\end{cases*}
	\hspace{0.5cm} \textup{and} \hspace{0.5cm} d = 4 - \sum_{i=1}^{4} k_{i}^{2} - \prod_{i=1}^{4} k_{i}, 
\end{equation*}
such that the projective closure $X \subset \Pbb^{3}_{\Qbb}$ of $U = \mathcal{U} \times_{\Zbb} \Qbb$ is smooth. Then we have
\begin{equation}
  \#\{ k = (k_{1},k_{2},k_{3},k_{4}) \in \Zbb^{4}, |k_{i}| \leq M\; \forall\; 1 \leq i \leq 4 :  \mathcal{U}(\textbf{\textup{A}}_{\Zbb}) \not= \emptyset \} \asymp M^{4}  
\end{equation}
and also
\begin{equation}
  \#\{ k = (k_{1},k_{2},k_{3},k_{4}) \in \Zbb^{4}, |k_{i}| \leq M\; \forall\; 1 \leq i \leq 4 : \mathcal{U}(\textbf{\textup{A}}_{\Zbb}) \not= \mathcal{U}(\textbf{\textup{A}}_{\Zbb})^{\textup{Br}_{1}} \not= \emptyset \} \asymp M^{4}  
\end{equation}
as $M \rightarrow +\infty$.
\end{theorem}

\begin{proof}
Denote by $\mathcal{S}$ the set in question and $\mathcal{M}$ is the set of $k$ such that $|k_{i}| \leq M$ for $1 \leq i \leq 4$. Hence, using the same notation as in Proposition 4.4, the density we are concerned with is given by

$$ \mu(\mathcal{S}) = \lim_{M \rightarrow \infty} \frac{\#(\mathcal{S} \cap \mathcal{M})}{\# \mathcal{M}} = \lim_{M \rightarrow \infty} \frac{\#(\mathcal{S} \cap \mathcal{M})}{M^{4}}. $$ 

We apply Proposition 4.4 with the polynomials $f_{i} = x_{i}(x_{i}^{2}-2)(x_{i}^{2}-4)/3 \in \Zbb[x'_{1},x'_{2},x'_{3},x'_{4}]$ for $1 \leq i \leq 4$ in which we write $x_{1} = 144x'_{1} + 127$, $x_{j} = 144x'_{j} + 5$ for $2 \leq j \leq 4$ (resp. $x_{i} = 144.p_{0}^{3} x'_{i} + r_{i}$ where the $x_{i}$ is actually the $k_{i}$ and $p_{0} \geq 11$ is the chosen prime in Lemma 4.2 and the residues $r_{i}$ satisfy the hypotheses of the lemma), and using the inclusion-exclusion principle we can compute that $c_{2} = 0$, $c_{3} = 0$, $$ c_{p,1}= C^{2}_{4}.3^2.p^{2} - (3.C^{3}_{4}.3^3.p + 3.3^{4}) + (4.3^{3}.p + 16.3^4) - C^{4}_{6}.3^{4} + C^{5}_{6}.3^{4} - 3^{4} = 27(2p^{2} - 8p + 9) $$ if $p > 3$ and $p \equiv \pm 3\,(\textup{mod}\,8)$ and $$ c_{p,2} = C^{2}_{4}.5^2.p^{2} - (3.C^{3}_{4}.5^3.p + 3.5^{4}) + (4.5^{3}.p + 16.5^4) - C^{4}_{6}.5^{4} + C^{5}_{6}.5^{4} - 5^{4} = 25(6p^{2} - 40p + 75) $$ if $p > 3$ and $p \equiv \pm 1\,(\textup{mod}\,8)$ \big(resp. except for $p_{0}$ with $p_{0} \equiv \pm 3$ mod $8$, since $k_{2} \equiv k_{3}^{2} - 2 \equiv k_{4}^{2} - 2$, $k_{3} \equiv k_{4} \not\equiv 0$ and $k_{2}^{2}-4 \not\equiv 0$ (mod $p_{0}$), we have $k_{i}(k_{i}^{2} - 2)(k_{i}^{2} - 4) \not\equiv 0$ mod $p_{0}$ for all $2 \leq i \leq 4$, and so $c_{p_{0}} = 0\big).$
Applying $(24)$ gives us a positive density in $(25)$:

$$ \mu(\bigcap_{i \not= j} \mathcal{R}_{i,j}) = \frac{1}{144^{4}} \prod_{\substack{p > 3 \\ p \equiv \pm 3\,(\textup{mod}\,8)}} \left(1 - \frac{c_{p,1}}{p^{4}} \right) \prod_{\substack{p > 3 \\ p \equiv \pm 1\,(\textup{mod}\,8)}} \left(1 - \frac{c_{p,2}}{p^{4}} \right), $$
and also in $(26)$:
$$ \mu(\bigcap_{i \not= j} \mathcal{R}_{i,j}) = \frac{1}{144^{4}.p_{0}^{12}} \prod_{\substack{p > 3, p \not= p_{0} \\ p \equiv \pm 3\,(\textup{mod}\,8)}} \left(1 - \frac{c_{p,1}}{p^{4}} \right) \prod_{\substack{p > 3 \\ p \equiv \pm 1\,(\textup{mod}\,8)}} \left(1 - \frac{c_{p,2}}{p^{4}} \right). $$

Finally, we need to consider the number of surfaces which are singular (see necessary and sufficient conditions given in Lemma 3.6). By Lemma 3.6, the total number of $k$ with $|k_{i}| \leq M$ for $1 \leq i \leq 4$ such that the surfaces are singular is just $O(M^{3})$ as $M \rightarrow \infty$, hence it is negligible. Therefore, we obtain that $\mu(\mathcal{S}) \geq \mu(\bigcap_{i \not= j} \mathcal{R}_{i,j}) > 0$.
\end{proof}

\begin{remark}
Continuing from the previous remark, it would be interesting if one can find a way to include Assumption 3.3 into the counting result, which would help us consider the Brauer--Manin set with respect to the whole Brauer group instead of only its algebraic part.
\end{remark}

\section{Further remarks}
In this section, we compare the results that we obtain in this paper with those in the previous papers studying Markoff surfaces, namely \cite{GS22}, \cite{LM20} and \cite{CTWX20}. 

First of all, for Markoff surfaces, we see from \cite{LM20} that given $|m| \leq M$ as $M \rightarrow +\infty$, the number of counterexamples to the integral Hasse principle which can be explained by the Brauer--Manin obstruction is $M^{1/2}/(\log M)^{1/2}$ asymptotically. This implies that \emph{almost all} Markoff surfaces with a nonempty set of local integral points have a nonempty Brauer--Manin set, and for the surfaces (relative character varieties) that we study in this chapter, a similar phenomenon is expected to occur (looking at the order of magnitude in the main counting result above). However, at present, we are not able to compute the number of these surfaces for which the Brauer--Manin obstruction can or cannot explain the integral Hasse principle in a similar way as in \cite{LM20} and \cite{CTWX20}. If we can solve (partly) this problem, it will be really significant and then we may see a bigger picture of the arithmetic of Markoff-type cubic surfaces.

\subsection{Markoff descent and reduction theory}
Recall that, in order to show that the integral Hasse principle fails in \cite{CTWX20}, the authors also make use of the fundamental set, or \textit{box}, in \cite{GS22} as a very useful tool to bound the set of integral points significantly and then use the Brauer group elements to finish the proofs. This mixed method has been the most effective way to prove counterexamples to the integral Hasse principle for Markoff surfaces which cannot be explained by the Brauer--Manin obstruction until now. However, in the case of Markoff-type cubic surfaces which we consider in this paper, the bounds do not prove themselves to be so effective when considered in a similar way. Indeed, let us recall the \textbf{Markoff descent} below for convenience; see \cite{GS22} and \cite{Wha20} for more details on the notation. Given $k = (k_{1},\dots,k_{n}) \in \Abb^{n}(\Cbb)$, we write
$$ \textup{H}(k) = \max\{1,|k_{1}|,\dots,|k_{n}|\}. $$

\emph{Surfaces of type $(1,1)$.} Let $\Sigma$ be a surface of type $(1,1)$. By Section 2, we have an identification of the moduli space $X_{k} = X_{k}(\Sigma)$ with the affine cubic algebraic surface in $\Abb^{3}_{x,y,z}$ given by the equation
$$ x^{2} + y^{2} + z^{2} - xyz - 2 = k. $$
The mapping class group $\Gamma = \Gamma(\Sigma)$ acts on $X_{k}$ via polynomial transformations. Up to finite index, it coincides with the group $\Sigma'$ of automorphisms of $X_{k}$ generated by the transpositions and even sign changes of coordinates as well as the Vieta involutions of the form $(x,y,z) \mapsto (x,y,xy - z)$.

The mapping class group dynamics on $X_{k}(\Rbb)$ was analyzed in detail by Goldman as discussed in \cite{Wha20}, and the work of Ghosh--Sarnak \cite{GS22} (Theorem 1.1) establishes a remarkable exact fundamental set for the action of $\Sigma'$ on the integral points $X_{k}(\Zbb)$ for admissible $k$. More generally, the results below establish Markoff descent for complex points. 

\begin{lemma} \textup{\cite[Lemma 4.2]{Wha20}} 
Let $\Sigma$ be a surface of type $(1,1)$. There is a constant $C > 0$ independent of $k \in \Cbb$ such that, given any $\rho \in X_{k}(\Sigma,\Cbb)$, there exists some $\gamma \in \Gamma$ such that $\gamma^{*}\rho = (x,y,z)$ satisfies
$$ \min\{|x|, |y|, |z|\} \leq C \cdot \textup{H}(k)^{1/3}. $$
\end{lemma}

\emph{Surfaces of type $(0,4)$.} Let $\Sigma$ be a surface of type $(0,4)$. By Section 2, we have an identification of the moduli space $X_{k} = X_{k}(\Sigma)$ with the affine cubic algebraic surface in $\Abb^{3}_{x,y,z}$ given by the equation
$$ x^{2} + y^{2} + z^{2} + xyz = ax + by + cz + d $$
with $a,b,c,d$ appropriately determined by $k$. The mapping class group $\Gamma = \Gamma(\Sigma)$ acts on $X_{k}$ via polynomial transformations. Let $\Sigma'$ be the group of automorphisms of $X_{k}$ generated by the Vieta involutions
\[
\begin{aligned}
\tau^{*}_{x} &: (x, y, z) \mapsto (a - yz - x, y, z),\\
\tau^{*}_{y} &: (x, y, z) \mapsto (x, b - xz - y, z),\\
\tau^{*}_{z} &: (x, y, z) \mapsto (x, y, c - xy - z).
\end{aligned}
\]
Two points $\rho, \rho' \in X_{k}(\Cbb)$ are $\Sigma'$-equivalent if and only if they are $\Sigma$-equivalent or $\rho$ is $\Sigma$-equivalent to all of $\tau^{*}_{x}\rho', \tau^{*}_{y}\rho'$, and $\tau^{*}_{z}\rho'$. 

\begin{lemma} \textup{\cite[Lemma 4.4]{Wha20}} 
Let $\Sigma$ be a surface of type $(0,4)$. There is a constant $C > 0$ independent of $k \in \Cbb^{4}$ such that, given any $\rho \in X_{k}(\Cbb)$, there exists some $\gamma \in k$ such that $\gamma^{*}\rho = (x, y, z)$ satisfies one of the following conditions:
\begin{enumerate}
	\item[(1)] $\min\{|x|, |y|, |z|\} \leq C$,
	\item[(2)] $|yz| \leq C \cdot \textup{H}(a)$,
	\item[(3)] $|xz| \leq C \cdot \textup{H}(b)$,
	\item[(4)] $|xy| \leq C \cdot \textup{H}(c)$,
	\item[(5)] $|xyz| \leq C \cdot \textup{H}(d)$.
\end{enumerate}
\end{lemma}

Clearly, the restrictions in this lemma are \emph{weaker} than those in the previous lemma, and so seems their effect. It would be interesting if one can find a way to apply these fundamental sets to produce some family of counterexamples to the integral Hasse principle which cannot be explained by the Brauer--Manin obstruction. A natural continuation from our work would be to find some sufficient hypothesis for $k$ such that the Brauer--Manin set of the general Markoff-type surface is nonempty, as inspired by \cite[Corollary 5.11]{LM20}, which we will discuss in some particular cases in the next part. Ultimately, similar to the case of Markoff surfaces, it is still reasonable for us to expect that the number of counterexamples which cannot be explained by the Brauer--Manin obstruction for these relative character varieties is asymptotically \emph{greater} than the number of those which can be explained by this obstruction. 

\begin{example}
Take $k_{1} = 127, k_{2} = 5, k_{3} = 5 + 144.5 = 725, k_{4} = 5 + 144.10 = 1445$, then $k = (k_{1},k_{2},k_{3},k_{4})$ satisfies Assumption A and we now have an explicit example of a Markoff-type cubic surface (where local integral points exist) given by the equation
$$ x^{2} + y^{2} + z^{2} + xyz = 1048260x + 187140y + 99300z - 667871675. $$
By using the box as discussed above, we will prove that this example gives no integral solution, i.e., it is a counterexample to the integral Hasse principle.

Indeed, from the proof of \cite[Lemma 4.4]{Wha20}, we can find that the constant $C$ independent of $k \in \Cbb^{4}$ may take the value $C = \textbf{48}$ in the condition $(1)$ and $C = \textbf{24}$ in all the other conditions $(2),(3),(4),(5)$ (thus we may choose \textbf{48} to be the desired constant in the statement of the Lemma). With this information, we can run a program on SageMath \cite{SJ05} to find integral points satisfying one of the restrictions $(2), (3), (4), (5)$, along with the help of Dario Alpern's website (Alpertron) \cite{Alp} to find integral points whose one coordinate satisfies the restriction $(1)$. More precisely, SageMath shows that there is no integral point in any of the boxes defined by $(2), (3), (4), (5)$, while Alpertron deals with conic equations after fixing one variable bounded in $(1)$ by transforming them into homogeneous quadratic equations and showing that some corresponding modular equations do not have solutions. As a result, we are able to prove that there is clearly no integral point in all those five cases, hence no integral point on the corresponding Markoff-type cubic surface.
\end{example}

\subsection{Some special cases of Markoff-type cubic surfaces}
Since $a = b = c = 0$ in $(2)$ implies that either three of the $k_{i}$ are equal to $0$ or $k_{1} = k_{2} = k_{3} = -k_{4}$ or any other permutation of $k_{i}$ (with $1 \leq i \leq 4$), the only cases that our Markoff-type cubic surfaces recover the original Markoff surfaces are given by equations of the form
$$ x^{2} + y^{2} + z^{2} + xyz = 4 - k_{0}^{2}, $$
or 
$$ x^{2} + y^{2} + z^{2} + xyz = (2-k_{0}^{2})^{2} $$
for some $k_{0} \in \Zbb$. For the former equation, it is an interesting Markoff surface which can be studied similarly as in previous work, with a remark that for any odd $k_{0}$ the equation will not be everywhere locally solvable as the right hand side is congruent to $3$ mod $4$ (see \cite{GS22}). For the latter equation, it always has the integral solutions $(\pm(2-k_{0}^{2}),0,0)$. Although these special cases are only of magnitude $M^{1/4}$ compared to the total number of cases that we consider, it is clear that our family of examples do not recover these particular surfaces. We will discuss it here in a more general situation.

Naturally, one can find different hypotheses from that of Assumption A to work for more general cases, since Assumption A exists for technical reasons as well as specific counting results. More precisely, we will now consider the special case when $k_{1} \not= k_{2} = k_{3} = k_{4}$ are integers such that the total field extension $[E:\Qbb] = 4$ and the set of local integral points on the Markoff-type cubic surface is nonempty. Interestingly, we are under the same framework as that of the general case in \cite[Section 3]{CTWX20}, to compute the algebraic Brauer group of the affine surface, which gives us 
$$ \textup{Br}_{1}\,U/\textup{Br}_{0}\,U \cong (\Zbb/2\Zbb)^{3} $$
with three generators
$$ \{ (x-2,k_{2}^{2}-4), (y-2,k_{2}^{2}-4), (z-2,k_{2}^{2}-4) \}. $$
Following \cite[Proposition 5.7 and Lemma 5.8]{LM20} or \cite[Lemmas 5.4 and 5.5]{CTWX20}, we have 
$$ \{ (x-2,k_{2}^{2}-4)_{2}, (y-2,k_{2}^{2}-4)_{2}, (z-2,k_{2}^{2}-4)_{2} \} = \{0, 1/2, 1/2\} $$
and 
$$ \{ (x-2,k_{2}^{2}-4)_{3}, (y-2,k_{2}^{2}-4)_{3}, (z-2,k_{2}^{2}-4)_{3} \} = \{ 0, 0, - \} $$
as multisets. As $k_{2}$ are odd, so $k_{2}^{2}-4 \equiv 5$ mod $8$, hence there is no way to describe this number by the form $3v^{2}$ for $v \in \Zbb$, which is favorable to give Hasse failures as in previous work. Due to the complexity of the integral values of the polynomial $X^{2}-4$ and their prime divisors, for now we can only deduce the vanishing of Brauer--Manin obstructions to the integral Hasse principle using a similar method as in \cite{LM20}.

\begin{proposition}
Let $k_{1} \not= k_{2} = k_{3} = k_{4}$ are integers satisfying the congruence conditions $k_{1} \equiv -k_{2} \equiv -5$ mod $16$ (resp. mod $9$), and the divisibility conditions: $(k_{1},k_{2}) = 1$, $(k_{1},k_{2}^{2}-4) = (k_{1}^{2}-4,k_{2}) = 1$, and $(k_{1}^{2}-2,k_{2}^{2}-2) = 1$, such that $[E:\Qbb] = 4$. Moreover, assume that there is a prime $p > 3$ such that $p$ divides $k_{2}^{2}-4$ to an odd order and $k_{1} \equiv -k_{2}$ mod $p$. Then there is no algebraic Brauer--Manin obstruction to the integral Hasse principle, but there is a Brauer--Manin obstruction to strong approximation for $\mathcal{U}$. 
\end{proposition}

\begin{proof}
Following similar arguments as in the proof of Proposition 4.1, we can prove that the set of local integral points is indeed nonempty. Note that under our assumption, $a = b = c \equiv 0$ mod $p$ and $d = (2-k_{2}^{2})^{2} \equiv 4$ mod $p$, so $ax + by + cz + d \equiv 4$ mod $p$ in the defining equation $(1)$ of $\mathcal{U}$. Now for $p > 5$, we let $\mathcal{B} = \langle (x-2,k_{2}^{2}-4), (y-2,k_{2}^{2}-4), (z-2,k_{2}^{2}-4) \rangle$ and follow the same arguments in the proof of \cite[Proposition 5.5]{LM20} to prove that the map 
$$ \mathcal{U}(\Zpbb) \rightarrow \textup{Hom}(\mathcal{B}, \Qbb/\Zbb), \hspace{1cm} \textbf{u} \mapsto (\beta \mapsto \textup{inv}_{p}\,\beta(\textbf{u})), $$
induced by the Brauer--Manin pairing, is surjective. Now we need to show that $\emptyset \not= \mathcal{U}(\textbf{\textup{A}}_{\Zbb})^{\textup{Br}_{1}} \not= \mathcal{U}(\textbf{\textup{A}}_{\Zbb})$. Let $\textbf{u} \in \mathcal{U}(\textbf{\textup{A}}_{\Zbb})$. Then by surjectivity, there exist $\textbf{u}_{p} \in \mathcal{U}(\Zbb_{p})$ such that $\textup{inv}_{p}\,\beta(\textbf{u}_{p}) = - \sum_{v \not= p} \textup{inv}_{v} \beta(\textbf{u})$ for all $\beta \in \mathcal{B}$ and $\textbf{u}_{p}' \in \mathcal{U}(\Zbb_{p})$ such that $\textup{inv}_{p}\,\beta(\textbf{u}_{p}') = 1/2 - \sum_{v \not= p} \textup{inv}_{v}\,\beta(\textbf{u})$ for some $\beta \in \mathcal{B}$. Then we obtain a point in $\mathcal{U}(\textbf{A}_{\Zbb})^{\textup{Br}_{1}}$ by replacing the $p$-part of $\textbf{u}$ by $\textbf{u}_{p}$, and a point in $\mathcal{U}(\textbf{A}_{\Zbb})$ but not in $\mathcal{U}(\textbf{A}_{\Zbb})^{\textup{Br}_{1}}$ by replacing the $p$-part of $\textbf{u}$ by $\textbf{u}_{p}'$, as required.

We are left with the case $p = 5$. Using \cite[Proposition 5.7]{LM20} and \cite[Lemma 5.5]{CTWX20}, we know that the image of the above map induced by the Brauer--Manin pairing contains all the nontrivial elements. Let $\textbf{u} \in \mathcal{U}(\textbf{\textup{A}}_{\Zbb})$. In fact, there always exists $\textbf{u}_{p}' \in \mathcal{U}(\Zbb_{p})$ such that $\textup{inv}_{p}\,\beta(\textbf{u}_{p}') = 1/2 - \sum_{v \not= p} \textup{inv}_{v}\,\beta(\textbf{u})$ for some $\beta \in \mathcal{B}$. Next, if $\sum_{v \not= p} \textup{inv}_{v} \beta'(\textbf{u}) \not= 0$ for some $\beta' \in \mathcal{B}$, there still exists $\textbf{u}_{p} \in \mathcal{U}(\Zbb_{p})$ such that $\textup{inv}_{p}\,\beta(\textbf{u}_{p}) = - \sum_{v \not= p} \textup{inv}_{v} \beta(\textbf{u})$ for all $\beta \in \mathcal{B}$. Now if $\sum_{v \not= p} \textup{inv}_{v} \beta(\textbf{u}) = 0$ for all $\beta \in \mathcal{B}$, we can consider another local point $\textbf{u}''$ whose $p$-parts ($p \not= 2$) are the same as those of $\textbf{u}$ and $2$-part $(x_{2}'',y_{2}'',z_{2}'')$ is a permutation of $\textbf{u}_{2} = (x_{2},y_{2},z_{2})$ (note that the Markoff-type cubic equation here is symmetric in $x,y,z$) such that their images under the local invariant map at $2$ are \emph{different} permutations of $(0,1/2,1/2)$, hence we will get $\sum_{v \not= p} \textup{inv}_{v} \beta(\textbf{u}) \not= 0$ for some $\beta \in \mathcal{B}$. The proof is now complete.
\end{proof}

Again, the hypothesis of Proposition 5.3 can be modified or generalized to give the same results, which we hope to achieve in possible future work. For now, let us give some concrete examples from Proposition 5.3 using the help of SageMath \cite{SJ05} and Alpertron \cite{Alp} which have been mentioned in Example 5.1.

\begin{example}
Consider $k_{1} = 144.7.2 - 5 = 2011$ and $k_{2}=k_{3}=k_{4}=5$. Then there is no integral point on the corresponding Markoff-type cubic surface. By Proposition 5.3 (with $p = 7$), we obtain a counterexample to the integral Hasse principle which cannot be explained by the \emph{algebraic} Brauer--Manin obstruction. 
\end{example}

\begin{example}
Consider $k_{1} = 144.7.3 - 5 = 3019$ and $k_{2}=k_{3}=k_{4}=5$. Then we find an integral point $(x,y,z) = (24,409,672)$ on the corresponding Markoff-type cubic surface. By Proposition 5.3 (with $p = 7$), we get a counterexample to strong approximation which can be explained by the Brauer--Manin obstruction.
\end{example}

Finally, we end with a counterexample to the integral Hasse principle for which it is unclear whether the (algebraic) Brauer--Manin obstruction exists or not.

\begin{example}
Consider $k_{1} = 127$ and $k_{2}=k_{3}=k_{4}=5$. Since $(127^2-4,5^2-4) = 3$, this example does not completely satisfy any of our previous assumptions. However, by using the programs on SageMath and Alpertron as discussed above, we find that there is indeed no integral point on the corresponding Markoff-type cubic surface.
\end{example}

\selectlanguage{english}
\printbibliography


\textsc{Sorbonne Université and Université Paris Cité, CNRS, IMJ-PRG, F-75005 Paris, France}\\
\textit{E-mail address}: \href{mailto:quang-duc.dao@imj-prg.fr}{\texttt{quang-duc.dao@imj-prg.fr}}
\\~\\
\indent \textsc{Institute of Mathematics, Vietnam Academy of Science and Technology, Cau Giay, 122300 Hanoi, Vietnam}\\
\textit{E-mail address}: \href{mailto:dqduc@math.ac.vn}{\texttt{dqduc@math.ac.vn}}

\end{document}